\documentclass[11pt,reqno]{amsart}
\usepackage[dvipsnames]{xcolor}
\usepackage[left=0.8in,top=1in,right=0.8in,bottom=1in]{geometry}
\usepackage{mathtools, amssymb, tikz, mathrsfs, cases}

\usepackage[T1]{fontenc}
\usepackage[colorlinks=true, linkcolor=orange, urlcolor=violet, citecolor=blue]{hyperref}
\usetikzlibrary{matrix,arrows}

\hypersetup{
 pdfauthor={P. Corvaja, A. Turchet and U. Zannier}}

\input xy 
\xyoption{all} %

\usepackage{lineno}

\usepackage{color}

\theoremstyle{plain}
\newtheorem{thm}{\protect\theoremname}[section]
 \newcommand\thmsname{\protect\theoremname}
 \theoremstyle{plain}

  \theoremstyle{remark}
  \newtheorem{rem}[thm]{\protect\remarkname}
  \theoremstyle{definition}
  \newtheorem*{example*}{\protect\examplename}
  \theoremstyle{definition}
  
  \theoremstyle{plain}
  \newtheorem{lem}[thm]{\protect\lemmaname}
  \theoremstyle{plain}
  \newtheorem{prop}[thm]{\protect\propositionname}
  \theoremstyle{plain}
  \newtheorem{cor}[thm]{\protect\corollaryname}
  \newtheorem*{thm*}{\protect\theoremname}
 
\makeatletter
  \theoremstyle{definition}
  \newtheorem{my@rem}[thm]{Remark}
  \renewenvironment{rem}{\begin{my@rem}}{\end{my@rem}}
\makeatother

\makeatother

  \providecommand{\examplename}{Example}
  \providecommand{\lemmaname}{Lemma}
  \providecommand{\propositionname}{Proposition}
  \providecommand{\remarkname}{Remark}
  \providecommand{\theoremname}{Theorem}
\providecommand{\theoremname}{Theorem}
 \providecommand{\corollaryname}{Corollary}

\newtheorem{theorem}{Theorem}[section]
\newtheorem{lemma}[theorem]{Lemma}
\newtheorem{proposition}[theorem]{Proposition}

\newtheorem{claim*}{Claim}
\newtheorem{theo}[theorem]{Theorem}

\newcounter{cor}
\newtheorem{thmy}{Theorem}

\newcounter{theo}
\newtheorem{thmyy}{Theorem}

\theoremstyle{definition}

\newcommand{\A}{{\mathbb A}}

\newcommand{\PP}{{\mathbb P}}
\newcommand{\C}{{\mathbb C}}

\newcommand{\Q}{{\mathbb Q}}
\newcommand{\R}{{\mathbb R}}

\newcommand{\Z}{{\mathbb Z}}

\newcommand{\calC}{{\mathcal C}}

\newcommand{\calS}{{\mathcal S}}

\newcommand{\calX}{{\mathcal X}}

\newcommand{\OO}{{\mathcal O}}

\DeclareMathOperator{\GL}{GL}

\numberwithin{equation}{section}
\numberwithin{table}{section}

\author{Pietro Corvaja}
\address[Corvaja]{Dipartimento di Scienze Matematiche, Informatiche e Fisiche\\
Universit\`a di Udine\\
Via delle Scienze, 206\\
33100 Udine --  Italy\\}
\email{pietro.corvaja@dimi.uniud.it}
\author{Amos Turchet}
\address[Turchet]{Dipartimento di Matematica e Fisica \\
Universit\`a di Roma Tre\\
Largo San L. Murialdo, 1\\
00146 Roma -- Italy}
\email{amos.turchet@uniroma3.it}

\author{Umberto Zannier}
\address[Zannier]{Scuola Normale Superiore\\
Piazza dei Cavalieri, 7\\
56126 Pisa -- Italy}
\email{umberto.zannier@sns.it}
\title{Rational distances from given points in the plane}
\date{\today}


\makeatletter
\@namedef{subjclassname@2020}{\textup{2020} Mathematics Subject Classification}
\makeatother


\begin{document}

 \begin{abstract}
    In this paper we study sets of points in the plane with rational distances from $r$  prescribed points  $  P_1,\ldots,P_r$.  A crucial case arises for  $r = 3$, where we provide  simple necessary and sufficient conditions for the density of this set in the real topology. We show in Theorem \ref{T.main} that these conditions can be checked effectively (via congruences), proving that a related class of $K3$ surfaces satisfies the local-global principle.
     In particular, these conditions are always satisfied   when $P_1,P_2,P_3$ are rational.
     This result completes and goes beyond the analysis  of Berry \cite{Berry}, who worked under stronger assumptions, not always fulfilled for instance in all the cases where $P_1,P_2,P_3$ are rational. 
    On the other hand, for $r\ge 4$, we show that points with rational distances  correspond to rational points in a surface of general type, hence conjecturally not Zariski dense. However, at the present, we lack methods to prove this, given the fact  that the surface is simply-connected, as we shall show. We give explicit proofs as well as describe in detail the geometry of the surfaces involved. In addition we discuss certain analogues for points with distances in certain ring of integers.
 \end{abstract}

\maketitle

\section{Introduction}

In this paper we continue the investigation started by the third author in \cite{Z} on sets of lattice points in the plane with integral distance to each point in a prescribed finite  set. We focus here mainly on points having {\it rational} distances rather than integral distances, though we shall also consider some aspects of integrality. We shall further  consider  variants of the problem where we allow the points and the distances to take value in different rings or fields, such as number fields and their ring of integers, e.g. the ring $\Z[i]$  of Gaussian integers.

The investigation  of subsets of the plane with integral distances is a fascinating historical   subject that appears already in a paper of Anning and Erd\H{o}s \cite{AE} and motivated the so-called Erd\H{o}s-Ulam problem that is still an open question; see \cite{ABT, Tao, Z} for discussions and further references. We also refer to Guy's collection \cite{Guy} for related problems, raised by several authors.

\smallskip

We denote by $d(P,Q)$ the Euclidean distance between the points $P$ and $Q$ in $\R^2$.
Let $P_1,\ldots,P_r$ be fixed  points on the Euclidean plane $\R^2$; we shall usually suppose they are in general position, i.e. pairwise distinct and three-by-three non aligned; however, we shall treat in detail also the case  $r=3$ and $P_1,P_2,P_3$ aligned. The main focus for the present article is the set of points with rational distances to $P_1,\ldots,P_r$. Not surprisingly, this set becomes sparser  when $r$ increases; conjecturally, it is not Zariski-dense whenever $r\geq 4$.

\smallskip 

The case $r= 2$ is easy; the set of points of the plane having  rational distances from two given points $P_1,P_2$ is dense in the real (hence also in the Zariski) topology.

\smallskip

  In the case $r=3$, T. G. Berry \cite{Berry} proved, under certain conditions, that the set is dense (in the real topology): he required that all square distances $d(P_i,P_j)^2$ are rational and at least one of the distances is rational. Completing his work, we provide an `if and only if' criterion for density:
  
  \begin{thm}\label{T.main}
  Let $P_1,P_2,P_3\in \R^2$ be three non-aligned points. The following are equivalent:
  
  \begin{itemize}
  
  \item[(i)] the set of points $Q\in\R^2$ having rational distance from $P_1,P_2,P_3$ is dense in $\R^2$;
  
  \item[(ii)] the quadratic form 
 \begin{equation}\label{E.quadr-form}
  (x,y)\mapsto \|x\overrightarrow{P_1P_3}+y\overrightarrow{P_2P_3}\|^2
  \end{equation}
  is defined over $\Q$ and represents a square.
  
  \item[(iii)] the square distances $d(P_i,P_j)^2$, $i,j=1,2,3$, are all rational and there exists at least one point $Q\in\R^2$ having rational distance from $P_1,P_2,P_3$;
  \end{itemize}
  \end{thm}
 
\medskip

Thanks to the Hasse-Minkowski theorem for quadratic forms, our condition $(ii)$ can always be checked by verifying the solvability of finitely many congruences. Also, it is plainly satisfied when $P_1,P_2,P_3$ are rational points, unlike Berry's one. 

The following statement provides a further condition equivalent to $(ii)$, and hence to $(i)$ and $(iii)$. We first recall the Iwasawa decomposition in $\GL_2^+(\R)$, the connected component of the identity of the group of invertible $2\times 2$ matrices with real coefficients: every  matrix $T\in\GL_2^+(\R)$ can be uniquely decomposed as a product
$$
T=A\cdot D\cdot U
$$
where $A\in\mathrm{O}(2)$ is orthogonal, $D$ is diagonal with positive entries on the diagonal and $U$ is upper triangular unipotent.
\smallskip

{\bf Addendum to Theorem \ref{T.main}}. Let $P_1,P_2,P_3$ be three non-aligned points on the plane. Choose an ordering such that the matrix $T=(\overrightarrow{P_1P_3} ,\overrightarrow{P_2P_3})$ has positive determinant. Then $(ii)$ is equivalent to the condition 

\begin{itemize}

\item[(iv)]  $D=\left(\begin{matrix} \sqrt{r} & 0 \\ 0&\sqrt{s}\end{matrix}\right)$ with $r,s$ positive rational numbers such that the quadratic form $rx^2+sy^2$ represents $1$ and $U$ has rational entries.
 
\end{itemize}
\smallskip

Let us call {\it admissible} a triangle $(P_1,P_2,P_3)$ for which condition (ii) (hence also condition (iv)) is satisfied. The family of admissible triangles   is acted on by a natural group of affine transformations. Of course, this set is invariant by all translations by vectors of the group $\R^2$.  We find it remarkable that it is also invariant by the affine  transformations whose differential lies in the group $\mathrm{GL}_2(\Q)$;  moreover, it is also invariant by the homotheties of a factor $\lambda$, where $\lambda^2$ sum of two squares in $\Q$.  It turns out that the group generated by these three kinds of transformations  is the maximal  group of affine transformations leaving invariant the set of admissible triangles; it can be characterised as the group of affine transformations 
\[
\binom{x}{y} \mapsto T\cdot \binom{x}{y}+ \binom{u}{v},
\]
where $T\in\GL_2(\R)$, $(u,v)\in\R^2$ and $T$ is a product $\lambda\cdot A$, with $A\in\GL_2(\Q)$ and $\lambda\in\R^*$ such that $\lambda^2$ is a sum of two rational squares. Let us call admissible such an affine transformation. 
 The action of the group of admissible transformations on the set of admissible triangles is not transitive, as shown by the following example:

\medskip

\noindent{\tt Example}: the triangles $(0,(1,0),(0,1))$ and $(0,(1,0),(0,\sqrt{3})$, where the symbol $0$ denotes the origin of the plane, are both admissible, but no admissible affine transformation sends one triangle into the other. The first triangle is admissible since it has rational vertices; the second one leads to the quadratic form $x^2+3y^2$ which represents a square. To see that the two triangles do not belong to the same orbit under the mentioned group of affine transformations, it suffices to check one by one the six affine transformations of the plane  sending the first triangle to the second one.  Since the group of admissible transformations acts transitively on the vertices of the first triangle, it suffices to check only one transformation sending the first triangle onto the second one, e.g. the affine transformation fixing the points $0,(1,0)$ and sending $(0,1)$ to $0,\sqrt{3})$. This transformation is not admissible.

\smallskip

One can look at the triangles $P_1,P_2,P_3$ satisfying the first part of condition $(ii)$, namely the rationality of the quadratic form \eqref{E.quadr-form}. Choosing for simplicity $P_3$ equal to the origin, the triple $P_1,P_2,P_3$ is identified   by a non-singular $2\times 2$ matrix $T\in\GL_2(\R)$, whose column are the vectors $\overrightarrow{P_1,P_3}, \overrightarrow{P_2,P_3}$. The rationality condition $(ii)$ amounts to the rationality of the symmetric matrix ${}^t T\cdot T$.  A set of matrices of this form is given by those of the form $A\cdot T_0$, with $T_0\in\GL_2(\Q)$ and $A$ in the orthogonal group $\mathrm{O}(2)$. 

In turn one could search for a complete characterization of such matrices (which will be given below): for instance, if the points are defined over a number field $K$ (say Galois over $\Q$ with Galois group $G$), the condition amounts to 
\begin{equation}\label{E.cocycle}
{}^t T\cdot T = {}^t T^\sigma\cdot T^\sigma
\end{equation}
for each $\sigma\in G$, i.e. $T\cdot {T^{\sigma}}^{-1}$ belongs to the orthogonal group $\mathrm{O}(2)$.  So the map $\sigma\mapsto  T^{\sigma} T^{-1}$ is a $1$-cocycle of $G$ with value in the orthogonal group. If the first cohomology group $\mathrm{H}^1(G,\mathrm{O}(2,{K}))$  would vanish, each such matrix would be of the above mentioned shape   $A T_0$ (i.e. the triangle $P_1,P_2,P_3$ would be isometric to a rational one).

However, the above example with the triangle $0,(1,0), (0,\sqrt{3})$ shows that this is not the case.

See below for a connection with this fact with the field of moduli for the surfaces $\calS_3$ associated to the triangles $P_1,P_2,P_3$.  
\smallskip

\smallskip

An important ingredient in the proof of Theorem \ref{T.main} will be the following result, of easier nature, about rationality of the  squares of the distances:

\begin{thm}\label{T.square-distances} 
Given three non-aligned points $P_1,P_2$ and $P_3$, the following are equivalent: 
\begin{itemize}

\item[(i)] the set of points $Q\in\R^2$  having rational squared-distance from $P_1,P_2,P_3$ is dense;

\item[(ii)] the quadratic form 
 \begin{equation} 
  (x,y)\mapsto \|x\overrightarrow{P_1P_3}+y\overrightarrow{P_2P_3}\|^2
  \end{equation}
  is defined over $\Q$.
  
  \item[(iii)] the square distances $d(P_i,P_j)^2$, $i,j=1,2,3$, are all rational.
  
\end{itemize}

 Moreover, under the above conditions, the set of points having rational squared-distance from $P_1,P_2,P_3$  is equal to the set of points of the form $P_3+x\overrightarrow{P_1P_3}+y\overrightarrow{P_2P_3}$, with $x,y\in\Q$.
\end{thm}

This result will be proved in section \ref{proof-T-square}.

\smallskip

Arguing as in Theorem \ref{T.main}, we can express (ii) above in terms of Iwasawa decomposition:

\smallskip

{\bf Addendum to Theorem \ref{T.square-distances}}. Setting $T=(\overrightarrow{P_1P_3},\overrightarrow{P_2P_3} )\in\GL_2(\R)$, and $T=A\cdot D\cdot U$ be its Iwasawa decomposition, the above conditions are equivalent to  

\begin{itemize}
\item[(iv)]  $D=\left(\begin{matrix} \sqrt{r} & 0 \\ 0&\sqrt{s}\end{matrix}\right)$ with $r,s$ positive rational numbers  and $U$ has rational entries.
\end{itemize}

\medskip

Going back to our original problem concerning rationality of the distances, the case $r=4$   becomes much harder: for instance  no point $Q$ is known  having rational distances from all the vertices of the unit square. We shall see that the set of points having rational distances from four given points in general position corresponds to the set of rational points on a certain algebraic surface of general type. According to the Bombieri-Lang conjecture, such rational points should be contained in a finite union of algebraic curves.

\smallskip

Note that if we replace the unit square by the unit pentagon (meaning a regular pentagon with sides of length $1$) our condition (ii) of Theorem \ref{T.square-distances} is not satisfied for any choice of three of the five vertices, hence even the set of points with rational square distances from the vertices is degenerate (actually it can be proven to be empty). The same holds for $n$-sided regular unit polygon for $n\geq 7$ (see \cite{n-gon} for a complete proof).
\medskip

To see the link between rational points on algebraic surfaces and the solutions to our problem, let us introduce the following notation.

\smallskip

Given $r\geq 3$ points $P_1,\ldots,P_r$ in general position, let $\calS_{P_1,\ldots,P_r}=\calS_r$ be the surface in $\A^r$ consisting in the set of elements $(z_1,\ldots,z_r)\in\A^r$ such that there exists a point $Q\in\A^2$ such that
$$
(d(P_i,Q))^2=z_i^2\qquad i=1,\ldots,r.
$$  
Clearly, such a point $Q$ is unique, so the surface $\calS_r$ is endowed with a projection $\calS_r\to \A^2$ sending $(z_1,\ldots,z_r)$ to the corresponding point $Q$.

\smallskip

Working over the complex number field, or even over the reals, the surface   $\calS_r$ could also be defined as the subsets of $\A^{2+r}$ formed by the points $(Q,z_1,\ldots,z_r)$ satisfying the above equation. However, our arithmetic problem consists in finding points with $(z_1,\ldots,z_n)\in\Q^r$, regardless of the rationality of the point $Q$. Hence our definition.

Our main results can be summarized as follows:

\begin{enumerate}

\item The surfaces $\calS_3$ have smooth models which are elliptic $K$3 surfaces. Under the necessary and sufficient condition $(ii)$ above, their rational points are dense in the real topology (proven in section \ref{sec:3pts}). 

\item  Smooth models for the surfaces $\calS_r$, for $r\geq 4$, are simply connected and of general type (proven in section \ref{sec:4pts}). According to the Bombieri-Lang conjecture, their set of rational points are degenerate, that is,  not dense in the Zariski topology.

We remark at once that no example is known to date of a simply connected surface of general type when one can prove this fact.

\item Concerning integral distances: the set of points in the plane with coordinates in a ring $R$ of $S$-integers in a number field, and having distances from two fixed points that belongs to  $R$,  is Zariski dense for instance when $R $ is a real or imaginary quadratic ring. This fits with the Vojta's conjectures (see \cite{BG}). 

\end{enumerate}

\smallskip

The condition that the quadratic form \eqref{E.quadr-form} is defined over $\Q$ will be proved to be equivalent to the condition that the surface $\calS_3$ is defined over $\Q$; hence its necessity for the Zariski-density of its rational points.

The condition that the quadratic form $(x,y)\mapsto \|x\overrightarrow{P_1P_3}+y\overrightarrow{P_2P_3}\|^2$ represents a square is equivalent to the existence of at least one point $Q$ having rational distance from $P_1,P_2,P_3$. 

It is worth noticing that varying $P_1,P_2,P_3$ among triples with mutual rational square distances one obtains a family of $K3$ surfaces with the following arithmetic property: whenever one such a surface contains a rational point outside  certain given curves ``at infinity'', since condition $(iii)$ of Theorem \ref{T.main} is satisfied this surface  contains a Zariski-dense set of rational points.  This phenomenon, which occurs for Brauer-Severi varieties and for certain Del Pezzo surfaces, is rather unusual in the case of $K3$ surfaces.  \smallskip

We note that the surfaces $\calS_r$ depend on the $r$ points $P_1,\ldots,P_r$ only up to the action of the group of similitudes\footnote{By similitude we mean the composition of a homothety with an isometry.} on the affine plane; hence over the complex number field $\C$, the surfaces  $\calS_r$  form an algebraic family of dimension  $2r-4$. In particular, the surfaces $\calS_3$ form a two-dimensional family (of Kummer $K$3 surfaces).\smallskip

We shall show (see the last Remark in section \ref{sec:4pts}) that the conjectured degeneracy the rational points of $\calS_r$, for $r \geq 4$, would answer the famous Erd\H{o}s-Ulam problem in the negative, i.e. every infinite  set of points having mutual rational distances would be contained in a line or a circle. 
In particular we obtain, a different proof of the result proved by Tao in \cite{Tao}, namely that the Bombieri-Lang conjecture implies a negative answer to the Erd\H{o}s-Ulam problem.

The conditional result hinges on the fact that the set of points under consideration corresponds to rational points in a simply connected surface of general type (see Proposition \ref{prop:SimpConn}). This is relevant since no example is known when the rational points can be proved (unconditionally) to be degenerate (over every number field) for such a surface of general type.  \smallskip

Concerning the case $r\leq 2$, while the original problem has a trivial solution (the set of points having rational distances from $P_1,P_2$ is always dense) an interesting problem arises by requiring that $Q$ is also a rational points. We are then led to studying the surface $ \calX_2$ defined as
$$
\calX_2=\{(Q,z_1,z_2)\in\A^2\times \A^2\, :\, d(P_i,Q)=z_i,\, i=1,2 \},
$$
and investigating the distribution of its rational points.  

We shall prove the following

\begin{theorem}
The surface $\calX_2$ is a Del Pezzo surface of degree $4$. Whenever $P_1,P_2$ are rational, it is rational over $\Q$, so in particular its set of rational points is Zariski-dense.
\end{theorem}

\smallskip

\texttt{Acknowledgments.} We thank Fabrizio Catanese for several discussions about the geometry and the topology of the surface $\calS_4$, and to provide us with the argument (and the reference) for Proposition \ref{prop:SimpConn}; thanks to Francesco Veneziano for pointing out reference \cite{n-gon}. We thank the anonymous referees for several comments that greatly improved our presentation, and in particular to point out the results of \cite{Berry}. AT is partially supported by the PRIN 2022 project 2022HPSNCR: Semiabelian varieties, Galois representations and related Diophantine problems and the PRIN 2020 project 2020KKWT53: Curves, Ricci flat Varieties and their Interactions and is a member of the INDAM group GNSAGA. PC is a member of the INdAM group GNSAGA and is supported by the Advanced grant ``Hyperbolicity in Diophantine Geometry''.

\medskip
\section{Preliminary geometrical considerations} \label{S.Intermezzo}
Following the notation introduced above, given points $P_1,\ldots,P_r$ on the plane $\A^2(\R)=\R^2$, we define the algebraic surfaces $\calX_r$ and $\calS_r$ by letting
$$
\calX_r:\quad \{(Q,z_1,\ldots,z_n)\in \A^{r+2}\,\: \, (d(Q,P_i))^2 = z_i^2\}\subset\A^{r+2}
$$
and, for $r\geq 3$,
$$
\calS_r:\quad \{(z_1,\ldots,z_n)\in \A^{r}\,\: \, \exists Q\in\Q^2, \, (d(Q,P_i))^2 = z_i^2\}.
$$
We shall frequently denote by the same symbols different projective or affine models of such surfaces. We note that for $r\geq 3$ and $P_1,\ldots,P_r$ in general position the surfaces $\calX_r$ and $\calS_r$ are isomorphic over $\C$.

\medskip
{\tt About the surface $\calX_1$}. Let us consider the points in the plane with rational distance to the origin. This leads to the Pythagorean equation $x^2 + y^2 = z^2$, which we view as the relation between the coordinates of a point $Q=(x,y) \in \A^2$ and its distance $z = d(O,Q)$ from the origin.  A rational solution corresponds to a rational point with rational distance from the origin. Such equation defines a cone in $\A^3$, hence a singular affine surface; it is the quadratic cover of $\A^2$ ramified over the pair of complex lines defined by the equation $x^2+y^2=0$. Of course, it is a rational surface, indeed also rational over $\Q$, and since the ancient times we know the rational parametrization given by formulas 
\[
x = 2\lambda p \qquad y = \lambda(p^2 - 1) \qquad z = \lambda(p^2 + 1).    
\]

However, in view of the study of the problem of points with rational distance from three given points, we would like to study also a different model. 

A non-singular projective model, which will be denoted by $\calX_1$, can be obtained as follows: noticing that the two aforementioned complex lines intersect at the origin, let us consider the blow-up of the projective plane over the origin, i.e. the point of coordinates $(x:y:t)=(0:0:1)$; the strict transform of the reducible conic $x^2+y^2=0$ is the union of two disjoint curves on such a surface, so a smooth divisor. The quadratic cover ramified over these curves is our surface $\calS_1$. 

Let us now look for explicit equations: the blow-up of the plane can be defined inside $\PP_2\times\PP_1$ as the set of solutions $(x:y:t),(\xi:\eta)$ to the linear equation
$$
(x+iy)\eta=(x-iy)\xi.
$$
Adjoining to the function field the square-root of the rational function $(x^2+y^2)/t^2$ amounts to extracting the square root of $\xi\eta$, or of $\xi/\eta$. Hence we can define our surface inside $\PP_2\times \PP_1\times\PP_1$ with coordinates $(x:y:z),(\xi:\eta),(z_0:z_1)$ by the equation above combined with the new equation
$$
\eta z_0^2=\xi z_1^2.
$$
Substituting in the above equation, i.e. projecting from $\PP_2\times\PP_1\times\PP_1\to \PP_2\times\PP_1$ by forgetting the $(\xi:\eta)$ projective line, we obtain the single equation in $\PP_2\times\PP_1$
$$
(x+iy) z_1^2=(x-iy) z_0^2.
$$
This is the equation of the surface $\calX_1$. Note that the projection to $(z_0:z_1)\in\PP_1$ endows $\calS_1$ with a structure of $\PP_1$-bundle over $\PP_1$: indeed the surface $\calS_1$ is the so-called Hirzebruch surface of degree two, abstractly defined as  
$\PP({\OO}_{\PP_1}\oplus {\OO}_{\PP_1}(2))$, i.e. the completion of the total space of a line bundle of degree two over $\PP_1$, obtained by adding a single point to each fiber. 
\smallskip

{\tt About the surface $\calX_2$} The above construction can be mimicked to study the surface $\mathcal{X}_2$. In this setting we are interested in points $Q \in \A^2$ with rational distance to two fixed points. Without loss of generality we can assume that one of the points is the origin $O = (0,0)$ and the other point is $P = (a,b)$, with $a,b \in \Q$. Then, as above, we consider a bi-quadratic cover of $\PP_2$ ramified over the pair of singular conics defined by the equations $x^2+y^2=0$ and $(x-a)^2+(y-b)^2=0$. Again, after blowing up the singular points of these conics, i.e. the   points $O$, $P$ in our original problem, we can replace the singular conics by a smooth divisor; the bi-quadratic covering ramified over such a divisor turns out to be a smooth surface, and indeed is a Del Pezzo surface of degree four (see the explicit computation in the next section).

\medskip

{\tt Integral points in $\calX_2$}. The distribution of integral points on such affine surfaces often represents a deep problem: it was implicitly considered by M. Davis in \cite{Davis} in his investigations on the Diophantine representation of the exponential functions (see also \cite{Cantone} for recent developments).

 In our setting, the  integral points on the plane having {\it integral} distances from the two given points can be viewed as the rational points on the (complete) surface $\calX_2$ which are integral with respect to a certain divisor, namely the pull-back of the divisor at infinity of $\PP_2$ in the covering ${\calX}_2\to\PP_2$. Denoting by $D$ this divisor, we note that $D$ lies in the anti-canonical class of ${\calS}_2$, so the open surface ${\calS}_2\setminus D$ is a so-called log--$K3$ surface. By a result of B. Hassett and Yu. Tschinkel, the integral points on surfaces obtained by removing an anti-canonical divisor from a Del Pezzo surface are {\it potentially} dense, i.e. the set of points with coordinates in a suitable ring of $S$-integers is dense (actually, in \cite{HTsch} this was proved under the assumption that the divisor at infinity is smooth, but S. Coccia removed this condition in the recent work \cite{Coccia2}). 
 
However, restricting to the ring $\Z$ of rational integers makes the situation different:  the third author, in \cite{Z} proved that under suitable simple necessary conditions  the Zariski-closure of the set of points with coordinates in $\Z$  of $\calX_2$ projects to the union of at most two lines and a non-empty finite union of hyperbolas.

\smallskip

On the other hand, the set of integral points having integral distances from three given points is expected to be degenerate over any number field; this would follow from Vojta's Conjecture, applied to an affine open set of a K3 surface.
\medskip

\section{Rational distances from two points}\label{sec:2points}

The goal of this section is proving by explicit computation the following result:

\begin{theorem}\label{T.2punti}
The set of rational points on the plane having rational distances from   two given rational points is Zariski-dense on the plane, and actually also dense in the real topology.
\end{theorem}
 
As before, without loss of generality, we can assume that one of the points is the origin $O$ and we will denote by $P=(a,b)$ the second point, for rational numbers $a,b \in \Q$. We are interested in studying the affine set $\calS_2$ of points $Q = (x,y) \in \A^2$ such that the distances $d(O,Q)$ and $d(P,Q)$ are rational, i.e. the set
\[
    \calX_2 = \calS_2(\Q) = \{ Q \in \A^2(\Q) : d(O,Q) \in \Q \text{ and } d(P,Q) \in \Q \}.
\]

We stress here that we do not ask that the points of $\calS_2$ have mutual rational distance, i.e. given $Q_1,Q_2 \in \calS_2(\Q)$ the distance $d(Q_1,Q_2)$ might  be, in general, an irrational number (so in particular $\calS_2$ is not a ``rational distance set'' in the sense of \cite{ABT,Tao}).

A point $Q \in \calX_2(\Q)$ corresponds to a rational solution to the system
    \begin{align}
        &x^2 + y^2 = z^2,  \label{eq:1_2pt} \\ 
        &(x - a)^2 + (y - b)^2 = (z - k)^2 \label{eq:2_2pt},
    \end{align}
where $z = d(O,Q)$, and $z -k = d(P,Q)$. The above presentation shows that the set $\calS_2$ is the set of $\Q$-rational points of the intersection of two quadrics in $\A^4$, namely the solutions of equations \eqref{eq:1_2pt} and \eqref{eq:2_2pt} respectively, where $x,y,z$ and $k$ are the coordinates of $\A^4$. The (projective completion of the) variety $\calS_2$ defined by these equations is a Del Pezzo surface of degree four, hence a geometrically rational surface. 
Following \cite{Z}, one can show that this variety is  rational over $\Q$, so its rational points are Zariski dense.

We can eliminate $z$ in the two equations as follows. First, one can subtract \eqref{eq:2_2pt} from \eqref{eq:1_2pt} to get
\begin{equation}\label{eq:4_2pt}
    2kz = 2 ax + 2 by - \delta,
\end{equation}
where $\delta = a^2 + b^2 - k^2$. Squaring \eqref{eq:4_2pt} and using \eqref{eq:1_2pt} we get 
\begin{equation}\label{eq:Ck}
    \calC_k: \quad (2ax + 2by - \delta)^2 = 4k^2(x^2 + y^2).
\end{equation}
The latter equation represents a family of affine conics $\calC \to B\simeq \mathbb{A}^1$ indexed by $k$. We note that  (at least for $k \neq 0$) rational points on a fiber $\calC_k$ correspond to points $Q$ with rational distances from $O$ and $P$ where the difference $d(O,Q) - d(P,Q) = k$ is fixed. In particular, if $\calC$ has a dense set of rational points then $\calS_2(\Q)$ will be Zariski-dense too; actually it will turn out that $\calS_2(\Q)$ is also dense in $\calS_2(\R)$ in the real topology. 
We observe that the projective completion of $\calC$ is a conic bundle $\overline{\calC} \to \overline{B}$ over a projective curve $\overline{B}\simeq \PP_1$. Moreover, to show that $\calC$ (and hence $\overline{\calC}$) has a dense set of rational points, it is sufficient to produce a rational section $B \to \calC$. 

To continue with our explicit computation, we start with a rational parametrization of \eqref{eq:1_2pt} with parameters $\lambda$ and $p$ such as
\begin{equation}\label{eq:thot}
    x = 2 \lambda p \qquad y = \lambda(p^2 -1) \qquad z = \lambda(p^2 + 1).
\end{equation}
Then we use \eqref{eq:2_2pt} with $x = a$ to get $k = z - y + b$. Finally we can express $\lambda = \frac{a}{2p}$ and $p = \frac{a}{k}-b$ thus expressing $x$ and $y$ (and $z$) in terms of $a,b$ and $k$ as wanted. This defines a rational section $B \to \calC$; since $B(\Q)$ is dense in $B(\R)$ in the real topology, and almost every  fiber of a rational point of $B$ has a dense (in the real topology) set of rational point,  we obtain that $\calC$ has a dense set of rational points.  Actually, the section shows that the surface is rational: for the section meets every conic in a rational point, and that conic can be parametrized rationally in terms of such rational point.

This proves Theorem \ref{T.2punti} that the set $\calX_2$, of rational points in the plane with rational distance from two fixed points $O$ and $P$ with rational coordinates is, in particular,   dense in $\A^2$, or, what amounts to the same, the set of rational points $\calX_2(\Q)$ is dense on the surface $\calX_2$: the density holds actually simultaneously with respect to any finite set of places; in other words, we can prescribe that our point with rational distance from both $P_1,P_2$ lies near to any prescribed point of the plane, and satisfies every prescribed congruence modulo any given integer. 

This is an explicit description parallel to the more algebro-geometric viewpoint that we described in section \ref{S.Intermezzo}, where we started from a geometrical description of the involved algebraic surface.

\medskip

\subsection{A third approach for distances from two points}\label{sec:2pt_conics}
 
The aim of this section is showing the geometrical meaning of the explicit calculations made in the previous section; we take the opportunity to reprove in an alternative and more synthetic way the density of the set of rational points having rational distances from two given rational points (i.e. the Zariski-density of $\calX_2(\Q)$).

\smallskip

{\tt Confocal families of conics}. Given two points $O,P$ on the plane $\A^2(\R)=\R^2$, let us  consider the algebraic family of conics (either ellipses or hyperbolae) having foci on $O$ and $P$. They form a $1$-parameter family, but not a linear system: every generic point of the plane belongs to two such conics. A confocal family can always be defined as the set of conics which are tangent to four fixed given lines; hence it is the dual family of a pencil of conics in the dual projective plane. 
 
 The four (complex) lines tangent to every conic of the family are the tangent drawn from the foci: those passing through $O$ are given by the equation $x^2+y^2=0$ (i.e. set of (complex) points with ``distance zero'' from $O$) and analogously for the other pair of lines, those passing through $P$).

Curiously, such considerations on families of confocal conics arose recently in a work by the first and third authors on finiteness results on elliptical billiards \cite{CZ-billiards}. 

The ellipses in a confocal family are defined as the set points $Q\in \R^2$ such that $d(O,Q)+d(P,Q)=k$, where $k$ is a given real number $>d(O,P)$. The hyperbolae in the family are defined by the relation $|d(Q,O)-d(Q,P)|=k$, where now $k$ is a real number with $0<k<d(O,P)$.

\smallskip

{\tt On the proof of Theorem \ref{T.2punti}}.  We note that the conics $\mathcal{C}_k$, defined by equation \eqref{eq:Ck} and used in the previous section  to prove our density statement, are precisely the conics having foci at $O$ and $P$. 

The fact that their rational points are solutions to our initial problem, i.e. they have rational distance from both $O$ and $P$ can be seen as follows: first of all the squared distances $(d(O,Q))^2, (d(P;Q))^2$ are rational numbers, and this happens for every rational point $Q$ on the plane, so is the difference 
$$
(d(O,Q))^2-(d(P,Q))^2=(d(O,Q)+d(P,Q))\cdot (d(O,Q)-d(P,Q)).
$$
Hence if one between the numbers $d(O,Q)+d(P,Q)$ and $d(O,Q)-d(P,Q)$ is rational, so is the other (provided neither is zero) and then both numbers $d(O,Q)$ and $d(P,Q)$ are rational. In geometric terms, if the point $Q$ belongs to either a hyperbola or an ellipse corresponding to a rational parameter $k$ as before, then $d(O,Q)$ and $d(P,Q)$ are rational numbers; also, if $Q$ belongs to a hyperbola of equation $|d(Q,O)-d(Q,P)|=k$ for some rational $k$, then it also belongs to a confocal ellipse of equation $d(Q,O)+d(Q,P)=k'$, where $k'$ is again a rational number.

We could conclude the proof of Theorem \ref{T.2punti} as follows: first we choose a single rational point $Q_0$   having rational distances from both $O$ and $P$, and note that  we can choose $Q_0$ outside the line segment joining $O$ and $P$. (This first step is easy; see e.g. the lemma below for a much stronger result.)  Then, consider the ellipse having foci at $O$ and $P$ passing through $Q_0$. This conic must contain infinitely many rational points, since it is smooth, defined over $\Q$ and contains at least one rational point by construction; for almost every point on this ellipse, so in particular for almost every rational point, there exists also a confocal hyperbola passing through it. By what we have just said, each rational point on each such hyperbola is a point with rational distance from both $O$ and $P$. This reproves Theorem \ref{T.2punti}.

\smallskip

In our proof in \S \ref{sec:2points} we followed a different pattern, implicitly based on the following lemma:

\begin{lemma}\label{L.retta}
Let $O\in \A^2(\Q)=\Q^2$ be a rational point and $l\subset \A^2$ a line defined over $\Q$ not passing through $O$. Then there exist infinitely many rational points $Q\in l(\Q)$ having rational distance from $O$.
\end{lemma}

\begin{proof}
We  choose coordinates in such a way that $O$ is the origin.  We may assume that $l$ is defined by $y=ax+b$, with rational $a,b$, where $b\neq 0$. The squared distance from $O$ of a point $(x,y)\in l$ is given by  $x^2+(ax+b)^2=(a^2+1)x^2+2abx+b^2$. For  $x=0$ such distance is rational, hence the conic   $z^2=(a^2+1)x^2+2abx+b^2$ has a rational point, and so it has infinitely many such points. (For instance set  $z=b+tx$; after dividing by $x\neq 0$, we get $t^2x+2bt=(a^2+1)x+2ab$ hence  $x=2b(a-t)(t^2-a^2-1)^{-1}$ is a parametrization in $t$). 

\end{proof}

To finish the proof of Theorem \ref{T.2punti}, we choose a point $Q_0\neq P$ having rational distance from $P$ and such that the line $l$ joining $P$ and $Q_0$ does not contain $Q$. Note that every rational point in $l$ has rational distance from $P$; by the above lemma, infinitely many rational points of $l$ have rational distance also from $O$. For each such point $Q$, consider the hyperbola having foci in $O,P$ passing through $Q$. Such hyperbola contains infinitely many rational points, each of them having rational distance from both $O$ and $P$. This is the strategy that we followed in the previous paragraph.

\smallskip

{\tt Yet another variant}. We note that we could have exploited Lemma \ref{L.retta}, working with an infinite family of lines passing through $P$, each of them possessing a rational point with rational distance from $P$ (for instance, take the unit circle centered on $P$, which has infinitely many rational point; each of them provides a suitable line); then, by Lemma \ref{L.retta}  each such line contains infinitely many rational points having rational distance from $O$, producing a dense set of rational points having rational distances from both $O$ and $P$. 

\medskip

{\tt A final geometrical remark}. The algebraic surface $\calX_2$ is obtained as a degree-4 cover of the $(x,y)$-plane, a Galois cover of type $(2,2)$, ramified over the two degenerate conics $x^2+y^2=0$ and $(x-a)^2+(y-b)^2=0$, $(a,b)$ being the coordinates of $P$. A generic conic on the plane lifts to a genus-one curve of the surface $\calS_2$. However, the conics of the 
  confocal family $\calC_k$ (see equation \eqref{eq:Ck}) considered in the above proof lift to rational curves on $\calS_2$; this is due to the fact that this family coincides precisely with the family of conics tangent to every component of the ramification divisor, which consists, as we said, in the union of four lines. Due to this tangency, the pre-image of each curve $\calC_k$ totally splits into four components, each one being a rational curve.

\medskip

\section{Squared Distances from Three Points}
In this section we begin the study of the set of points $Q \in \R^2$ that have rational distance from three fixed points, proving Theorem \ref{T.square-distances}.\medskip

{\tt Notation}.  Without loss of generality, we can assume that one of the three fixed points is the origin $O$. We will denote by $P=(a,b)$ and $P'=(a',b')$ the second and third point. 
For $Q, Q'\in \R^2$ we shall denote by $\|Q\|$ the usual euclidean length, i.e. $\|Q\| = d(Q,0)$  and by $Q.Q'$, or $(Q.Q')$, the associated scalar product.
Moreover, we will denote the three distances as follows $d_0 = \|Q\| = d(O,Q)$, $d = \|Q - P\| = d(P,Q)$ and $d' = \|Q - P'\| = d(P',Q)$.\medskip

In order to identify necessary and sufficient conditions for the density of the set of points with rational distances from $O,P$ and $P'$, we consider the set of points $Q\in\R^2$ such that the {\it squared distances} $d_0^2 = \|Q\|^2, d^2 = \|Q-P\|^2, d'^2 =\|Q-P'\|^2$ are all rational. 
 
In particular we suppose neither that $P,P'$, or $Q$, are rational points, nor that the distances between these points are rational. The output of this analysis, which is simple but apparently missing from previous papers in the topic, will lead to results holding essentially without restrictions. 
 
\medskip
 
We shall study conditions which ensure that the set  
  \begin{equation}\label{new:J}
  J:=\{Q\in\R^2:d_0^2,d^2,d'^2\in\Q\}.
  \end{equation}
is Zariski-dense in the plane $\R^2$. These are necessary conditions to hold in order for the set of points with rational distance from $O,P$ and $P'$ to be Zariski dense. 
 
 \medskip
 
\begin{rem}
It will follow from the analysis below that it suffices to assume that $J$ is not contained in a finite union of conics, in order to achieve the same conclusions. 
\end{rem}

 \medskip
 
 Setting $Q:=(x,y)$ we have that a point $Q \in J$ satisfies
\begin{align}
    &x^2 + y^2 = d_0^2  \label{new:eq:1} \\ 
    &(x - a)^2 + (y - b)^2 = d^2 \label{new:eq:2} \\
    &(x - a')^2 + (y - b')^2 = d'^2, \label{new:eq:3}
\end{align}

and, by construction of $J$, we presently assume that the right hand sides belong to $\Q$.

By subtraction, setting $g:=d_0^2-d^2$, $g'=d_0^2-d'^2$,   we get
\begin{equation}\label{new:line} 
2ax+2by=g+|P|^2,\qquad 2a'x+2b'y=g'+|P'|^2.
\end{equation}

\subsection{The collinear case} Suppose first that $O,P,P'$ are collinear, so $P'=q P$ for a $q \in \R^*\setminus\{1\}$.  Then 
\[g'= qg +(q\|P\|^2-\|P'\|^2)=q g+q(1-q)\|P\|^2.\]
 Hence $(g,g')$ lies on a certain fixed  line, not parallel to any of the axes.  Since $P\neq O$,  if $J$ is Zariski-dense the set of such $(g,g')$ is Zariski-dense on the line.  Hence the line is a rational line, which amounts to  $q\in\Q$, $|P|^2\in\Q$.  

Furthermore we can show that these conditions are also sufficient.

\begin{prop}  \label{new:Jdense_coll} Suppose that $O,P,P'$ are collinear: $P'=q P$, with $q\neq 0,1$. Then $J$ is Zariski-dense if and only if $p:=|P|^2\in \Q$ and $q\in\Q$. Under these assumptions, if $R = (-b,a)$ and $Q=tP+uR$, where $t=(P.Q)/p$ and $u=(R.Q)/p$, the necessary and sufficient conditions for $Q$ to be in $J$ are that $t\in\Q$ and $u^2\in\Q$.
\end{prop}

\begin{proof}
   The above discussion already proved the necessity.
   
   Conversely, let us assume that $q,\|P\|^2\in\Q$ (where $P'=qP$). First of all if $d_0^2\in \Q$ and if the first of \ref{new:line} holds for a $g\in\Q$, then the $g'$ obtained by the second of \eqref{new:line} will be rational as well, and $(x,y)$ will belong to $J$.  Hence we need only study \eqref{new:eq:1} and \eqref{new:eq:2} in rationals $d_0^2,d^2$.  In turn, this is like imposing that both $x^2+y^2=d_0^2$ and $ax+by=:\sigma$ are rational,  where we are assuming that $a^2+b^2$ is rational.  Now, for given $(d_0^2,\sigma)\in\Q^2$ we can solve for $x,y$, so the set $J$ is indeed Zariski-dense.  Specifically, we obtain the equations for $x,y$ given by 
\begin{equation*}
|P|^2x^2-2a\sigma x+\sigma^2-b^2d_0^2=0,\qquad |P|^2y^2-2b\sigma y+\sigma^2-a^2d_0^2=0.
\end{equation*}
Finally, if we let $R=(-b,a)$, we have that $\|R\|^2=\|P\|^2$ is rational and $(P.R)=0$.  Given $Q=tP+uR$ for $t,u\in\R$, we have that $|Q|^2$ and $(P.Q)$ have both to be rational. In this case, we have $|Q|^2=(t^2+u^2)|P|^2$ and $(P.Q)=t\|P\|^2$, thus obtaining the desired result.
\end{proof}

\begin{rem}
If we require that for instance the points in $J\cap\Q^2$ are still Zariski-dense then we find that $P,P'$ must be rational points as well; and this condition is also  sufficient.
\end{rem}

\subsection{ The non-collinear case. Proof of Theorem \ref{T.square-distances}}\label{proof-T-square} Suppose now that $0,P,P'$ are not collinear, so $\Delta:=ab'-a'b\neq 0$. 

Now, if $J$ is Zariski-dense we deduce from \eqref{new:line}  that the set $\Gamma\subset \Q^2$ of corresponding pairs $(g,g')$ is also Zariski-dense.  

Solving  the equations \eqref{new:line} for $x,y$ and inserting the values so obtained  in \eqref{new:eq:1} we obtain a condition  of the form $\Phi(g,g')\in\Q$, where $\Phi$ is a certain nonzero quadratic polynomial depending only on $P,P'$.  More precisely we have \begin{equation}\label{new:xy}
x=A(g+\|P\|^2)+B(g'+\|P'\|^2),  \qquad y=C(g+\|P\|^2)+D(g'+\|P'\|^2), 
\end{equation}
where
\begin{equation}\label{new:ABCD}
2\Delta A=b', \quad 2\Delta B=b, \quad 2\Delta C=-a', \quad 2\Delta D=-a.
\end{equation}

Setting ${\bf v}:=(A,C)$, ${\bf w}:=(B,D)$, we then have 
\begin{equation}\label{new:phi}
 \Phi(u,v)=\|{\bf v}\|^2u^2+2({\bf v}.{\bf w})  uv+\|{\bf w}\|^2v^2+Eu+Fv+G,
\end{equation}
where
\begin{equation*}
E=2\|{\bf v}\|^2\|P\|^2+2({\bf v}.{\bf w})\|P'\|^2, \quad  F=2\|{\bf w}\|^2\|P'\|^2+ 2({\bf v}.{\bf w})\|P\|^2,\quad  G=\|{\bf v}\|^2\|P\|^2+\|{\bf w}\|^2\|P'\|^2.
\end{equation*}

Now, if the vectors $(g^2,gg',g'^2,g,g',1)\in\Q^6$ lie in a proper vector subspace of $\Q^6$ (for varying $Q$ in $J$) then a fixed equation $\Psi(g,g')=0$ would hold, $\Psi$ being some nonzero quadratic polynomial, and, by \eqref{new:line}, $J$ would not be Zariski-dense.  

Therefore the vectors span $\Q^6$, and we deduce that $\Phi$ must have  rational coefficients. In the first place   this implies that the  three quantities $\|P\|^2/\Delta^2, (P.P')/\Delta^2, \|P'\|^2/\Delta^2$ (which are essentially  the coefficients of the homogeneous quadratic part of $\Phi$) are rational.  Now, let us look at $E,F$, which may be seen as linear forms in $\|P\|^2,\|P'\|^2$ with rational coefficients. The matrix of these linear forms up to a factor $2$  has rows $(\|{\bf v}\|^2|, {\bf v}.{\bf w})$ and $(\|{\bf w}\|^2,{\bf v}.{\bf w})$, which is the product of the matrix with rows ${\bf v},{\bf w}$  with its transpose $({\bf v}^t, {\bf w}^t)$. But this product matrix is rational, and it is also nonsingular, the  determinant of  $({\bf v}^t, {\bf w}^t)$ being $(4\Delta)^{-1}$. Hence $|P|^2,|P'|^2$ are rational, whence $\Delta^2$ is rational, and therefore $(P.P')$ is rational as well. 

\medskip

Note that these conditions amount simply to $P,P'\in J$; in fact, this is equivalent to $\|P\|^2,\|P'\|^2,\|P-P'\|^2$ to be all rational, and  the claim follows  in view of the equality $\|P-P'\|^2=\|P\|^2-2(P.P')+\|P'\|^2$.

\begin{rem}\label{rmk:Delta}
We note in passing that the rationality of $\|P\|^2,(P.P'),\|P'\|^2$ implies itself that of $\Delta^2$. In fact is $\alpha$ is the angle $POP'$ we have $(P.P')=\|P\|\|P'\|\cos\alpha$, $\Delta^2=\|P\|^2\|P'\|^2\sin^2\alpha=\|P\|^2\|P'\|^2-(P.P')^2$.
\end{rem}
 
Conversely, suppose that $\|P\|^2,(P.P'),\|P'\|^2\in\Q$. Then by Remark \ref{rmk:Delta} necessarily $\Delta^2\in\Q$ and the above argument can be reversed, giving a Zariski-dense set of points $Q$ with rational squared-distance from $O,P,P'$. 
 
 \medskip
 
 We now summarize these conclusions and describe in a simple way the set $J$.
 
 \begin{prop} \label{new:Jdense} Suppose that $0,P,P'$ are not collinear. If $J$ is Zariski-dense, then $\|P\|^2,(P.P'),\|P'\|^2$ are rational, or equivalently  $P,P'\in J$,  and this condition is also sufficient. Also,  if this holds we have $J=\Q P+\Q P'$, and actually $J$ is an orthogonal sum $J=\Q P+\Q R$ for a suitable $R\in J$ with $(P.R)=0$, $\|R\|^2\in\Q$.  Finally, the squared-distance between any two points in $J$ is rational. 
 \end{prop}

\begin{proof} Let $L:=\Q P+\Q P'$, a vector space of dimension $2$ over $\Q$, contained in $\R^2$. Since $\|P\|^2,(P.P')\in\Q$, we may find a nonzero vector $R\in L$ with $(P.R)=0$: it suffices to put $R:=P'-\frac{(P.P')}{\|P\|^2}P$. Then $L=\Q P+\Q R$. Now, if $Q\in L$, we have $Q=tP+uR$ with rational $t,u$ and then $|Q|^2=t^2|P|^2+u^2|R|^2\in \Q$, and similarly $|Q-S|^2\in\Q$ for every $S\in L$.  Also, if $Q\in J$, namely if $\|Q\|^2, \|Q-P\|^2,\|Q-P'\|^2\in\Q$, we have that $(Q.P),(Q.P')\in\Q$, whence $t,u\in\Q$, so $Q\in L$ as wanted. The last assertion also follows.
\end{proof} 

\begin{cor} If $J$ is Zariski-dense and contains two non collinear  points in $\Q^2$, then $J=\Q^2$.
\end{cor}

\begin{proof} Indeed, under these assumptions $P,P'$ cannot  be collinear with the origin, and we may apply the proposition. Then $J$ has a $\Q$-basis formed of rational points, whence the conclusion.
\end{proof}

\medskip

{\it Proof of the addendum to Theorem~\ref{T.square-distances}.} Let $P_1,P_2,P_3$ be a triangle satisfying condition (ii) of Theorem~\ref{T.square-distances}. Consider the vectors $\overrightarrow{P_1P_3}, \overrightarrow{P_2P_3}$ and suppose they form an oriented basis of $\R^2$. Consider the quadratic form $(x,y)\mapsto \|x\overrightarrow{P_1P_3}+y\overrightarrow{P_2P_3}\|^2$. We want to prove that it is defined over $\Q$ if and only if condition (iv) of the Addendum is satisfied. 

Since the euclidean quadratic form is invariant under the orthogonal group $\mathrm{O}(2)$, we can act with a matrix $A\in\mathrm{O}(2)$ in such a way that the vector $\overrightarrow{P_1P_3}$ takes the form ${}^t(a,0)$, where $a>0$, without changing the quadratic form $(x,y)\mapsto \|x\overrightarrow{P_1P_3}+y\overrightarrow{P_2P_3}\|^2$. The second vector $\overrightarrow{P_2P_3}$ will be of the form ${}^t(u,b)$ with $b>0$. The quadratic form reads
\[
a^2 x^2+2au xy + (u^2+b^2)y^2.
\]
Its rationality implies that $a$ is of the form $a=\sqrt{r}$, for some rational $r>0$ and $u=\sqrt{r}\cdot \xi$ for some $\xi\in\Q$; then $b^2\in\Q$, so we can write $b=\sqrt{s}$ for some rational number $s>0$; so the new matrix $(\overrightarrow{P_1P_3}, \overrightarrow{P_2P_3})$ takes the form $\left(\begin{matrix}\sqrt{r} & \sqrt{r}\xi\\ 0 & \sqrt{s}\end{matrix}\right)$. After multiplying this matrix by the inverse of the diagonal matrix $D=\left(\begin{matrix}\sqrt{r}&0\\ 0&\sqrt{s}\end{matrix}\right)$ we obtain a unipotent upper triangular matrix whose upper-right entry is the rational number $\xi$ as wanted.

\medskip

\section{Rational distances from three points}\label{sec:3pts}

Using the analysis from the previous section we consider now the set 
\begin{equation*}
J^*=\{Q\in\R^2: d_0:=\|Q\|\in\Q, d:=\|Q-P\|\in\Q,d':=\|Q-P'\|\in\Q\}.
\end{equation*}
As before, the points $P,P'$ are distinct points in the plane, which are also distinct from the origin $O$.

Of course $J^*\subset J$ (see Definition \ref{new:J}), and in fact {\it a priori} $J^*$ is a much smaller set. While the density of $J$ is a necessary condition for the density of $J^*$ we will show that if $J$ is Zariski-dense then $J^*$ too is Zariski-dense, provided certain very simple necessary conditions hold. In particular this will extend the main theorem of \cite{Berry} by providing the optimal set of necessary conditions that guarantee the density of points with rational distances from the given points $O,P$ and $P'$.

\subsection{The collinear case} Let us start again with the collinear case.  We put $p:=\|P\|^2$, $P'=qP$, and we assume that $p,q\in\Q^*$ which, as we have seen, amounts to $J$ being Zariski-dense. 

As above, we set $P=(a,b)$, $R=(-b,a)$ so $\|R\|^2=p$ and $(P.R)=0$. Then we have verified that $J$ consists of the points $tP+uR$ such that $t,u^2\in\Q$.   Now for a point $Q=tP+uR$ we have
$$
\|Q\|^2=(t^2+u^2)p,\qquad \|Q-P\|^2=((t-1)^2+u^2)p,\qquad \|Q-P'\|^2=((t-q)^2+u^2)p.
$$

 Hence by multiplying by $p$   we see that $Q$ lies in $J^*$ if and only if $t^2+u^2, (t-1)^2+u^2,(t-q)^2+u^2$ are of the shape $p$ times a rational square. Putting $v:=u^2\in \Q$ we arrive at the system
 $$
 t^2+v=pz^2,\qquad (t-1)^2+v=p(z+k)^2,\qquad (t-q)^2+v=p(z+k')^2,
 $$
 where $t,v,z,k,k'\in\Q$. By subtraction, on putting $v:=pz^2-t^2$,  this amounts to solving
$$
-2t+1=pk^2+2pzk,\qquad -2qt+q^2=pk'^2+2pzk',\qquad t,z,k,k'\in\Q.
$$
 In turn,  given $k,z$ we may define $t$ by the equation on the left, and we remain with the single equation
 $$
 pqk^2+2pqzk+q^2-q=pk'^2+2pzk',\qquad z,k,k'\in\Q.
 $$
This defines a family of conics parametrized by $z$. Since we may solve rationally for $z$, this defines a rational surface over $\Q$, whence in particular its rational points are Zariski-dense (and satisfy weak approximation).

\subsection{The non-collinear case} We now turn to the non-collinear case, which is less evident and involves  a condition which is not always verified. In fact, this is not surprising since the rationality in question amounts to study rational points on a surface which is not rational: indeed, it is a Kummer surface (see Theorem~\ref{T.Kummer}).

So let now $P,P'\in\R^2$ be non-collinear. 

In the previous section we have found a simple necessary and sufficient condition for $J$ to be Zariski-dense, given by  Proposition \ref{new:Jdense} therein. Using that result we may restrict our attention to points $Q$ in the vector space $J=\Q P+\Q R$, where $(R.P)=0$ and $r:=\|R\|^2\in\Q$. It follows from the proof of Proposition~\ref{new:Jdense} that we may take \[r=\|P'-\frac{(P.P')}{\|P\|^2}P\|^2=\|P'\|^2-(\frac{(P.P')}{\|P\|})^2=\frac{\Delta^2}{\|P\|^2}\].

Setting
\begin{equation}\label{new:pp's}
p:=\|P\|^2, \qquad   p':=\|P'\|^2,\qquad   s:=(P.P'),
\end{equation}
the quadratic form $\|Q\|^2$ with respect to the basis $P,R$ of $J$ takes the  shape 
\begin{equation*}
\Phi(x,y)=px^2+ry^2,
\end{equation*}
 if $x,y\in\Q$ are the coordinates of $Q$ in the said basis. As before, we may take $r=p'-(s^2/p)$.

\medskip

Note that  $|\cdot |$ may attain a nonzero rational value on $\Q P+\Q P'$ only if the quadratic form $\Phi$ represents a square over $\Q$, or (by definition) if the so-called Hilbert symbol $(p,r)=1$.  (Note also  that $(p,r)=(p,p\Delta^2)=(p,-\Delta^2)$.) This also amounts to the fact that the ternary form $\Phi(x,y)-z^2$ is isotropic (i.e. represents $0$) over $\Q$. 
We shall prove that this condition is sufficient. 

\begin{rem}
Note that {\it a priori}, given the condition,  it does not seem obvious even that there {\it exists} some point $Q$ with rational distance from $O,P,P'$, the condition stating only  {\it the existence of a point  $Q\in \Q P+\Q P'$ with rational distance from $O$}. Note also that the condition is obvious when for instance $P,P'$ are rational points, in which case $J=\Q^2$. See equation \eqref{new:existence} for explicit formulae. giving rational solutions.
\end{rem}

\begin{thm} \label{new:thm} Suppose that $P,P'\in\R^2$ are two non collinear points. If $\|P\|^2,(P.P'),\|P'\|^2\in\Q$, the set  $J^*$ is Zariski-dense if and only if there exists a  point $\Q P+\Q P'$ with nonzero  rational distance from the origin. On the other hand, if the opening assumption does not hold, then already $J$ is not Zariski-dense.
\end{thm}

To prove the theorem it will be convenient to work inside $\Q^2$ equipped with the form $\Phi$, which corresponds to work in $J$ with the basis $P,R$. 

\medskip

 In the said  basis $P$ has coordinates $(1,0)$ whereas we may choose $R$ as in the proof of Proposition~\ref{new:Jdense}, so that $P'$ has coordinates $(c,1)$ in this basis, for some $c\in\Q$. We shall work throughout with these coordinates.  However now the quadratic form $x^2+y^2$ is replaced by the new one. 
 
 Note that with this notation we have $s=(P.P')=(cP+R.P)=c\|P\|^2=cp$.

Let us write down the relevant equations for a point $Q=(x,y)\in J^*$. 
\begin{align}
    &px^2 + ry^2 = z^2  \label{eq:1} \\ 
    &p(x - 1)^2 + ry ^2 = (z-k)^2 \label{eq:2} \\
    &p(x - c)^2 + r(y - 1)^2 = (z-k')^2, \label{eq:3}
\end{align}

Similarly to the above, on subtracting the second equation from the first, we obtain a family of affine conics with parameter $k$ with fibers given by 

\begin{equation}\label{eq:C_k}
    \calC_k: \quad 2px - \delta  = 2kz,\qquad \delta=p-k^2.
\end{equation}

Using the first and third equation we also obtain a second family, with parameter $k'$:
\begin{equation}
    \label{eq:C'k'}
    \calC_{k'}': \quad 2pcx + 2 ry - \delta'  = 2k'z,\qquad \delta' = pc^2+r- k'^2.
\end{equation}

Let us write  down explicitly one possible cubic fibration. Setting $\lambda=y/x$ we obtain $y=\lambda x$, $z=\sqrt{p+r\lambda^2} x$ (for some choice of the square root)  and so the equations for $\calC_k, \calC_{k'}'$ become
\begin{equation}\label{eq:sistema}
\left\{ \begin{matrix}
2px-p &=& -k^2+2k\sqrt{p+r\lambda^2}\, x \\
(2pc+2r\lambda)x-pc^2-r &=&-k'^2+2k'\sqrt{p+r\lambda^2}\, x
\end{matrix}\right.
\end{equation}
 which correspond to  the compositum of two quadratic covers of the $x$-line, each one ramified  over two points. Such a curve has genus $1$.  Note that we can choose $\lambda$ to be rational and such that $p+r\lambda^2$ is a rational square in infinitely many ways. With  such choices the  coordinate $z=\sqrt{p+r\lambda^2}\, x$ also becomes rational; moreover both curves $\calC_k$ and $\calC_{k'}'$ will have a rational point above $x=\infty$, so that the genus-one curve defined by the system of equations \eqref{eq:sistema}  will have a marked rational point (actually four rational points), turning that curve into an elliptic curve (with rational $2$-torsion). 
 
 \medskip
 
 Let us derive a cubic curve out of these equations. 
 
 Write $\mu^2=p+r\lambda^2$. This is a conic, denoted $S$, which by assumption has a rational point, hence can be parametrized over $\Q$.\footnote{We have used affine coordinates, but of course we may refer to points at infinity on the conic, given by $(1:\pm\sqrt r:0)$ in homogeneous coordinates.}
  In terms  of these coordinates  the equations may be rewritten as 
 \begin{equation*}
\left\{ \begin{matrix}
2(p-k\mu )x  &=& p-k^2\\
2(pc+r\lambda-k'\mu)x &=& pc^2+r-k'^2.
\end{matrix}\right.
\end{equation*}
 
 Eliminating $x$ we obtain a family of  cubic equations, depending on the point $(\lambda,\mu)$ on our conic $S$:
 
 \begin{equation}\label{new:cubic}
 (p-k^2)(pc+r\lambda-\mu k')=(pc^2+r-k'^2)(p-\mu k).
 \end{equation}
 
 Note that if for rational $(\lambda,\mu)$ in our conic we have a rational point $(k,k')$ in the corresponding cubic, then, {\it provided} $p\neq \mu k$ and $ pc+r\lambda\neq \mu k'$, we can use the previous equations so as to obtain rational values for $x,y,z$, hence a rational point $Q$ with the sought properties.

\subsubsection{Some sections for the cubic family}

We also abbreviate $p':=pc^2+r=|P'|^2$, $\eta:=pc+r\lambda$. 

\medskip

{\bf A homogeneous equation for the cubic family}:  In homogeneous coordinates $(K:K':H)$ the above equation reads 
\begin{equation}\label{E.equazione-superficie}
 \mathcal{C}:\quad   (p H^2-K^2)(\eta H-\mu K')=(p'H^2-K'^2)(pH-\mu K).
 \end{equation}

It is the equation of a  cubic curve depending on a parameter which is a point of the  conic $S$ of equation $\mu^2=p+r\lambda^2$. The cubic is generically smooth because it is birationally equivalent to  the former curve, which  has genus $1$. (One can also check this by differentiation: one finds that possible singular points would have $k=k'$, which is easy to exclude. See also Lemma \ref{T.liscio} below.)

\smallskip

The affine coordinates $k,k'$ are obtained from the homogeneous ones by $k=K/H,\, k'=K'/H$.

\smallskip

{\bf Four rational sections}. There are  four rational sections, namely
$$
A:=(1:0:0),\, B:=(0:1:0),\, C:=(1:1:0),\, N:=(p:\eta:\mu).
$$
\footnote{Here $A,B,C$ are independent of $(\lambda,\mu)$ but we view them as  points with coordinates which are constant functions on $S$.}  
We note that, as is easy to check, $N$ does not coincide with any of $A,B,C$, no matter the point $(\lambda,\mu)$ on the conic (not even at infinity): this follows easily from the positivity of $|\cdot |^2$ on $\R^2\setminus\{O\}$. Note also that two of the four factors defining the cubic vanish at $N$.

Further, it is easily checked that these points are smooth on $\mathcal C$ for every point on the conic. 

\smallskip


It will be useful to note that
\begin{equation}\label{E.p^2}
\begin{matrix}
\mu^2 p-p^2 &= &(p+r\lambda^2)p-p^2 &=& pr\lambda^2 \\
\mu^2 p'-\eta^2 &= & (p+r\lambda^2)(pc^2+r)-(pc+r\lambda)^2 &=& pr(c\lambda-1)^2 \\
\end{matrix}
\end{equation}

\bigskip

\begin{lem}\label{Lemma1}
The tangents drawn from the  points  $A,B$  intersect on the fourth point. 
\end{lem}

\begin{proof} Let us prove that the tangent to the point $A=(1:0:0)$ passes through the fourth point.
De-homogenizing with $H=1$ the equation becomes
$$
(p-k^2)(\eta-\mu k')-(p'-k'^2)(p-\mu k)=0.
$$
The line joining $(1:0:0)$ with $(p:\eta:\mu)$ has equation $K'=H\eta/\mu$, i.e. $k'=\eta/\mu$ in affine coordinates. Substituting we obtain
$$
(p'-\eta^2/\mu^2)(p-\mu k)=0
$$
which has the only solution $k=p/\mu$. Since a line has a unique point at infinity, this means that there is a double solution at infinity on the cubic.

The verification for the point $(0:1:0)$ is symmetric.
\end{proof}

\begin{cor}
Taking the point $N$ as origin for the group law on the cubic curve, supposed to be smooth, we have identically on $S$
$$
A+C=B,\qquad B+C=A
$$
and so
$$
2C=0.
$$
 \end{cor}

This follows immediately from the linear equivalences provided by the lemma, together with the fact that $A,B,C$ are collinear points on the cubic: $A+B+C\sim 2A+N\sim 2B+N$.

\medskip

Let us now consider the point $(1:1:0)$. We have the following lemma.

\begin{lem}\label{Lemma-tangente}
Suppose $\mu\neq 0$. The line $l_C$ joining $C=(1:1:0)$ to $N=(p:\eta:\mu)$ is tangent to the cubic at the point $(1:1:0)$ if and only if $\lambda$ is such that $\eta=p$. It is  tangent at $N$ if and only if $c\lambda-1=\pm \lambda$. 
If $p=p'$ then $\eta=p$ implies $c\lambda-1=-\lambda$, and under such a specialization the cubic is reducible and the line $l_C$ is one of its components. If $p\neq p'$, for no specialization of $\lambda$ can the curve contain the line $l_C$ (except if $\mu=0$, i.e. $\lambda^2=-p/r$).   
\end{lem} 
 
\begin{proof}
The line joining the points $(1:1:0)$ to the point $(p:\eta:\mu)$ is parametrized as
$$
(p+t:\eta+t:\mu),\qquad t\in\mathbb{P}_1.
$$
Substituting into \eqref{E.equazione-superficie} we get a cubic equation with a factor $\mu t$. Indeed
the equation   simplifies to
$$
\mu t\cdot (p\mu^2 -p^2-2p t - p' \mu^2+\eta^2 +2\eta t)=0.
$$
 If $\mu=0$ the cubic is degenerate and contains the whole line $H=0$, a case which we are disregarding. The value $t=0$ yields the point $N$. The value $t=\infty$ yields $C$, whereas a third solution is obtained from the polynomial in brackets, which is 
 $
  2(\eta -p)t+ (p\mu ^2-p^2)-(p'\mu^2-\eta^2).
 $
 By \eqref{E.p^2} this equals 
 $$
 2(\eta -p)t+pr(\lambda^2-(c\lambda-1)^2).
 $$

Then, if $(\eta-p,\lambda^2-(c\lambda-1)^2)\neq (0,0)$ it gives  
$$
t=pr\frac{(c\lambda-1)^2-\lambda^2}{2(\eta-p)}.
$$
Then $t$ vanishes if and only if $c\lambda-1=\pm \lambda$ and is infinite if and only if $\eta=p$. In this second case the line is tangent to the point at infinity, proving the first assertion of the lemma.

Let us now analyze the cases when it happens that for some $\lambda$ both $\eta-p$ and $\lambda^2-(c\lambda-1)^2$ vanish. 

We can exclude that $\eta-p=0=c\lambda-1-\lambda$ because in this case $c\neq 1$ and we obtain that $p(c-1)^2+r=0$ which is impossible since $p,r>0$. On the contrary, it can happen that $\eta-p=0=(c+1)\lambda-1$. Indeed this amounts to $(c+1)\lambda=1$ and  to  $p(c^2-1)+r=0$, which in turn amounts to  $p=p'$ (and can never happen if $c=-1$).

In this case the specialization $\lambda\mapsto 1/(c+1)$ which makes $\eta$ be equal to $p$ leads to a reducible cubic, containing the line $l_C$.
\end{proof}

In the generic case (that is, for generic values of $\lambda$)  the line $l_C$ contains another point beyond $C,N$; since $div(l_C)\sim A+B+C$, such  point is  equal to the sum $A+B$   (with respect to the group law with origin $N$). The computations above exhibit this point:
\begin{equation}\label{new:existence}
A+B=(2p(\eta-p)+pr(c\lambda-1)^2-pr\lambda^2: 2\eta(\eta -p)+pr(c\lambda-1)^2-pr\lambda^2 : 2\mu (\eta-p)).
\end{equation}
 
\begin{rem}
 Letting $\lambda$ vary correspondingly to rational points on our conic, this section already provides an infinite set of rational solutions for our problem. Indeed, we find that $k=p/\mu$ or $k'=\eta /\mu$ only if $t=0$, which, as we have seen,  happens only if $(c\pm 1)\lambda=1$.
\end{rem}

Next, we want to prove that the sections found so far are not all torsion, with the purpose  to use the group law on the elliptic curves of the family in order to produce a Zariski-dense set of rational points.

In principle there are several methods for this, for instance one could specialize at some rational points; however there are parameters which complicate the situation. A suitable method is to explore points of the conic $S$ where a given section becomes torsion of low order: if this occurs at a point of good reduction, then either the section is torsion of {\it that} order, or it is non-torsion.

 With this in mind, 
 we look for  a specialization of $\lambda$ for which the elliptic curve has good reduction,  and an extra relation between $A,B,C$ holds; namely, we want that   $A:=(1:0:0)$ and $ B:=(0:1:0)$ become torsion of order $2$ with respect to the `moving' point $N:=(p:\eta:\mu)$. This will occur when the point $N$ is a flex and all the tangents drawn from $A,B,C$  contain $N$. In view of the lemma just proved, this happens if and only if $p=\eta$. However, we must check that this is a good specialization, namely that the specialized cubic has genus $1$.

\begin{lem}\label{T.liscio}
For no (complex) value of $\lambda$ the corresponding cubic curve is singular at the point $N=(p:\eta:\mu)$. 
The unique tangent to the cubic at $N$ meets the point $A=(1:0:0)$ if and only if $\lambda=0$ and meets $B$ if and only if $c\lambda=1$.
In these two cases the cubic becomes reducible and the tangent at $N$ is one of its components.
\end{lem}

\begin{proof}
Working in affine coordinates, i.e. putting  $H=1$, the equation of the cubic curve reads
$$
f(k,k'):=(p-k^2)(\eta-\mu k')-(p'-k'^2)(p -\mu k)=0,
$$
and the point $N$ in question is the point $(p/\mu, \eta/\mu)$. 
We have
$$
\dfrac{\partial f}{\partial k}= -2k (\eta-\mu k') +\mu(p'-k'^2), \qquad \dfrac{\partial f}{\partial k'}=-\mu (p-k^2)+2k'(p-\mu k);
$$
specializing at the point $(p/\mu, \eta/\mu)$ we obtain for the gradient the expression:
$$
\nabla f\left(\frac{p}{\mu},\frac{\eta}{\mu} \right) = \frac{1}{\mu}( \mu^2p' -\eta^2,-\mu^2p +p^2).
$$
As already observed  we have
$$
\mu^2 p -p^2=pr\lambda^2,\qquad \mu^2p' -\eta^2=pr(c\lambda-1)^2,
$$
arriving at
\begin{equation*}\label{E.gradiente}
\mu\nabla f \left(\frac{p}{\mu},\frac{\eta}{\mu} \right) =pr((c\lambda-1)^2, -\lambda^2).
\end{equation*}
Note that the vector in the right hand side is never zero. 
  This  proves that $N$ is a smooth point.

This formula also proves that  the (unique) tangent line has the direction of the vector $((c\lambda-1)^2, -\lambda^2)$, so its  point  at infinity is $A=(1:0:0)$ if and only if $\lambda$ vanishes and is $B$ if and only if $c\lambda-1'$ vanishes. 
 
 Finally, in these two cases such a tangent is a bitangent of the cubic by Lemma \ref{Lemma-tangente}, so it must be a component. 
\end{proof}

\subsubsection{\bf Good reduction} First of all, as we already observed, $\mu=0$ is a point of bad reduction, since for that value the cubic degenerates in the union of a conic and the line $H=0$. This is the only case in which the curve contains the line $H=0$. In this case, the second component is a smooth conic. (In fact, otherwise we have $p\eta=pp'$ and $\mu=0$, leading to $r\lambda^2=-p$ and $rp+(pc(c-1)+r)^2=0$, which is impossible since $p,r>0$.)

\medskip

Let us classify the other cases of reducibility.  Let $l_A$ (resp. $l_B, l_C$) denote the lines joining $A$ (resp. $B,C$) to the point $N$ and $l$ denote the line containing $A,B,C$ (i.e. the line $H=0$). Note that $l$ intersects the cubic precisely at  these three distinct points for every specialization of $\lambda$ with $\mu\neq 0$. Hence these points remain smooth for every specialization (outside the case $\mu=0$ when $l$ is a component of $\mathcal{C}$). By Lemma \ref{T.liscio} the point $N$ too remains smooth for every specialization. 

\smallskip

Suppose now that for a specialization of $\lambda$ the cubic becomes reducible. Let us first exclude it is the union of three lines.
We have already seen that $l$ is a component of $\mathcal{C}$ if and only if $\mu=0$, and in this case the second component is a smooth conic. So we can suppose that $l$ is not a component of $\mathcal{C}$.
Then $A,B,C$ are smooth points, as we observed; then, since the tangents from $A$ and $B$ must be unique and contain $N$, two of the three components would be $l_A,l_B$, making $N$ become a singular point, contrary to Lemma \ref{T.liscio}.

\medskip

Suppose now that $\mathcal{C}$ specializes to the union of a smooth conic and a line. Since we are excluding $\mu=0$, such a line cannot be $l$. Hence two among $A,B,C$ belong to the smooth conic and the other one to the line. 

If the line contains $A$, then it is $l_A$ and the smooth conic must contain $B,C$, not $N$ (which already belongs to $l_A$ and is smooth). This is the case $\lambda=0$.

Symmetrically, if the line contains $B$, then it is $l_B$ and $A,C$  belong to the smooth conic. This occurs if and only if $c\lambda-1=0$.

If the line contains $C$, let us show that it must be $l_C$, so by Lemma \ref{Lemma-tangente}  $p=\eta$ and $\lambda=-c\lambda+1$ and $p=p'$. Indeed, otherwise $A$ and $N$ would belong to the smooth conic. However $l_A$ would be tangent to that conic at $A$ and intersect that conic also at $N$.

Finally, the line cannot contain only the point $N$ among the four considered points, since the three aligned points $A,B,C$ cannot lie on a smooth conic.

\medskip

Finally, reducibility sometimes arises: for $\mu=0$ (the curve contains $l$),  for $\lambda=0$ (the curve contains $l_A$), for $\eta=0$ (the curve contains $l_B$) and, in the case $p=p'$, also  for $\lambda= - c\lambda+1$ (the curve contains $l_C$).

Recall that we are particularly interested in the specialization $p=\eta$. We now prove:

\begin{lem}\label{L.triangolo}
If  $p, p', s$ are pairwise distinct,  the specialization $p=\eta$, i.e. $\lambda=p(1-c)/r$, leaves the cubic curve irreducible.
\end{lem}

\begin{proof}
Suppose that $O,P,P'$ satisfy the hypotheses of the lemma.
By what we said, when the curve $\mathcal{C}$ becomes reducible, it contains one of the lines $l,l_A,l_B,l_C$.

 Since $p\neq p'$ by Lemma  \ref{Lemma-tangente} the curve specialized at $p=\eta$  does not contain $l_C$. 
 
 Suppose that it contains $l_A$ under the specialization $p=\eta$; by Lemma \ref{T.liscio} this means that $\lambda=0$, so $c=1$, whence (in the above notation) $P'=P+R$ where $(R.P)=0$. This implies $s=(P.P')=|P|^2=p$, contrary to assumptions.

 The argument for $B$ is symmetric: if the cubic contains $l_B$ (for $p=\eta$) then we derive from Lemma \ref{T.liscio} that $c\lambda=1$. This, together with $p=\eta$ yields $p'=pc^2+r=pc=s$.
 
 So far we have excluded that the curve contains $l_C,l_A,l_B$. 
 
 Let us exclude that it may contain $l$; otherwise $\mu=0$, so $\lambda^2=-p/r$. This is not compatible with $p=\eta$, which corresponds to a real value of $\lambda$.
 \end{proof}

 We now have a general geometrical lemma.
 
\begin{lem}\label{L.geometrico}
Let $\mathcal{C}$ be an irreducible cubic curve, $Q\in\mathcal{C}$ any point. Suppose that there exist three distinct smooth points $A,B,C\in\mathcal{C}\setminus\{Q\}$ such that the tangents to the curve at these points  all contain the point $Q$.  
 Then $\mathcal{C}$ is smooth. If moreover $A,B,C$ are collinear, then   $Q$ is a flex.
\end{lem} 
 
\begin{proof}
The projection from $Q$ provides a rational map $\mathcal{C}\to \PP_1$ of degree $2$ ramified (at least) at the points $A,B,C$. Then by Hurwitz genus formula the geometric genus of $\mathcal{C}$ is (at least) $1$, so it coincides with the arithmetic genus and the curve is smooth.

This proves the first assertion of the lemma.

Suppose now $A,B,C$ are collinear and consider the three lines $l_A,l_B,l_C$ joining the three points $A,B,C$ to $Q$, which by assumption are tangent to $\mathcal{C}$ at $A,B,C$ respectively. Denote also by $l_A,l_B,l_C$ three corresponding linear forms providing their equations. Let now $l$ be the line containing $A,B,C$ and again use the same symbol for a corresponding linear form. Finally, let $l^*$ denote both the tangent to the curve at $Q$ and a corresponding linear form.

The rational function $\displaystyle{\frac{l^2\cdot l^*}{l_A\cdot l_B\cdot l_C}}$, viewed as a function on the curve $\mathcal{C}$,  is regular and non vanishing at $A,B,C$; if $Q$ were not a flex, it would have a simple pole at the point $Q$ and no other pole, resulting in a  function of exact degree $1$, which is impossible on any smooth cubic curve.
\end{proof}

\smallskip

\medskip

{\bf The sections are not torsion}. As above, take for the origin the point  $N=(p:\eta:\mu)$ and consider the three sections provided by the (smooth) points $A:=(1:0:0),\, B:=(0:1:0),\, C:=(1:1:0)$.

We want to specialize $\lambda$ so that $p=\eta$; in that case the line $l_C$ will be tangent at $C$ and the values of the three sections will satisfy $A+B=C$ for the group law, since the linear equivalence $A+B+C\sim 2C+N$ then holds.

\medskip

\begin{prop} If   $p, p', s$ are pairwise distinct,   the two  sections $A$ and $B$ are not both torsion.
\end{prop}

Note that we always have  $A-B=C$, which is torsion, hence  that the two sections are not both torsion amounts   that none of them is torsion.  

\begin{proof}
First note that they are not identically torsion of order $2$: indeed if any of $A,B$  were identically torsion of order $2$, then $N$ would be identically a flex, since the tangents at $A,B$ pass through $N$,  so we would have $2A+N\sim 3N$ or $2B+N\sim 3N$, and the tangent at $N$ would not meet the curve elsewhere. That this is not the case can be proved by explicit computation, but we can also argue as follows. 

Suppose by contradiction that  both $A$ and $B$ are identically $2$-torsion; since $A+C=B$ we would have also identically $C=A+B$; however, the tangent drawn from $C$ does not always meet the origin $N$, which (in view of calculations above) prevents $C$ to be equal to $A+B$.   

\smallskip

Coming back to our main issue, let us now specialize $\lambda$ to the unique point such that $p=\eta$, i.e. $\lambda=p(1-c)/r$, so that by Lemma \ref{Lemma-tangente} the tangent at $C$ passes through $N$. 

Since we are assuming  the hypotheses of the previous Lemma \ref{L.triangolo}, the curve is irreducible. Then Lemma \ref{L.geometrico} asserts it is smooth, i.e. this value of $\lambda$ is a specialization of good reduction.

\smallskip

On the other hand, under this good specialization the sections $A,B$ satisfy the extra relation $A+B=C$ (which provides also $2A=2B=O$).

By general (elementary)  theory of reduction of elliptic schemes,  this proves that $A,B$ are not identically torsion, for otherwise the torsion order would be preserved by reduction, contrary to the opening assertion.
\end{proof}

\medskip

In view of this result, in order to achieve our goal (i.e. to find a non-torsion rational section), it suffices to work  on the assumption that  two among $p=|P|^2$, $p'=|P'|^2$, $s=(P.P')$ are equal. 

However note that  our goal remains unchanged if we permute the three points $O,P,P'$, and the same holds for our assumptions. Indeed, changing  the origin  of $\R^2$ amounts to add a same vector to $O,P,P'$ and this leaves invariant the space $\Q (P-O)+\Q (P'-O)$, which can be defined as the $\Q$-vector subspace of $\R^2$ made up of linear combinations of the three points, with zero sum of coordinates, which is independent of the origin.

Suppose then that the triangle $O,P,P'$ is not equilateral. If not all three sides have equal lengths, then   we may choose the origin so that $p\neq p'$. Also, we may assume that $p=s$.  Then $P'-P$ is orthogonal to $P$. and  $|P'-P|^2=p'-p$. If this is different from $p$, then we may move the origin in $P$ and obtain a triangle with two  distinct orthogonal sides, falling into the `good' case  already discussed. Therefore we may assume that $p'-p=p$, i.e. that  $p'=2p$. 

This corresponds to the fact that the triangle $OPP'$ is isosceles and right-angled  in $P$.   
Let us first treat this case.

\medskip

{\bf The case of the rectangular isosceles triangle}. In this case the strategy will be the same as above, but we shall use another specialization. Note that in the present situation we have $c=1$, $P'=P+R$, $r=p$, and $p'=2p$.

We take  $\lambda=1-c\lambda=1-\lambda$, i.e. $\lambda=1/2$, so that, , by Lemma \ref{Lemma-tangente} the line $l_C$ joining $C$ to $N$ becomes tangent  to the cubic at $N$. Since then $A+B+C\sim C+2N$, this   implies $A+B=0$, which is an extra relation between $A,B$. 

We have then to prove that this specialization leads to a smooth curve.

Note that in the present situation we have $\eta=pc+r\lambda=3p/2$, $\mu^2=p(1+\lambda^2)=5p/4$.

 Our cubic becomes (in affine coordinates)
 $$
 (p-k^2)(\frac{3p}{2}-\mu k')=(2p-k'^2)(p-\mu k).
 $$
 Note that it suffices to exclude singularities at finite points since first of all, the curve is smooth at infinity, since the line $H=0$ is not a component of the curve and it intersects the curve at three distinct points.

 We write the equation as $XY=ZW$. Differentiating with respect to $k$ and $k'$ we see that a  point is singular if and only if the two conditions $2kY=\mu Z$ and $\mu X=2k'W$ hold. 
 
 Suppose first that $W=0$, so $k=p/\mu$; then $k^2=p^2/\mu^2=4p/5\neq p$.  So $X\neq 0$ and we  the second condition does not hold. Hence $W\neq 0$. Multiplying by $W$ the first condition we get $2kWY=\mu ZW=\mu XY$. 
 
 If $Y=0$, then $Z=0$ as well and we find $k'=3p/(2\mu)$ and $k'^2=2p$, hence $9p^2/(4\mu^2)=2p$, which in turn yields $9p/5=2p$, which is impossible. 
 
 Hence $Y\neq 0$ and $2kW=\mu X$. Since $W\neq 0$ this yields $k=k'$, and we remain with the conditions
 $$
 (p-k^2)(\frac{3p}{2}-\mu k)=(2p-k^2)(p-\mu k),\qquad 2k(p-\mu k)=\mu(p-k^2).
 $$
 They respectively yield  $pk^2-2p\mu k+p^2=0$ and  $\mu k^2-2p k+\mu p=0$.  The discriminants up to a factor $4$ are $p^2\mu^2-p^3=p^3/4$ and $p^2-p\mu^2=-p^2/4$, hence one equation has two real roots, the other one two non-real conjugate roots. So there is no common solution whatever $p$.

 \bigskip
 
 We are left with the case of an equilateral triangle. Now $p=p'$, $s=p/2$, $c=1/2$, $r=5p/4$.

 {\bf The case of the equilateral triangle}. In this case  we choose the specialization $\lambda=c\lambda-1=(\lambda/2)-1$, i.e. $\lambda=-2$,  and prove it leads to a smooth cubic.  We have $\mu^2=p+4r=6p$.
 
 \smallskip

Again, we can work with the affine model, which takes the shape
$$
(p-k^2)(-2p-\mu k')=(3p/2 -k'^2)(p-\mu k).
$$
In order to compute possible singularities, we argue as before, writing this as $XY=ZW$. The two derivatives yield the equations $2kY=\mu Z$ and $\mu X=2k'W$. 

Again, suppose that $W=0$, then $X=0$, hence $p^2=\mu^2k^2=6pk^2=6 p^2$, which does not hold. Multiplying by $W$ the first condition and using the cubic equation we get 
$2kYW=\mu XY$. If $Y=0$ then $Z=0$, hence $(2p)^2=(\mu k')^2=6pk'^2=9p^2$, which once more does not hold. Hence we can divide out by $Y$ and find  $2kW=\mu X$, whereas also $\mu X=2k'W$ holds. Since $W\neq 0$ we get again $k=k'$, leading to the system
$$
(k^2-p)(2p+\mu k)=(3p/2-k^2)(p-\mu k),\qquad 2k(p-\mu k)=\mu(p-k^2),
$$
which may be rewritten as
$$
6pk^2+p\mu k-7p^2=0,\qquad \mu k^2-2pk+p\mu=0.
$$
The discriminants are respectively $p^2\mu^2+164p^3>0$ and $4p^2-4p\mu^2=-18p^2<0$. As before there cannot be common solutions, which concludes the verification.

\bigskip

{\bf End of the proof of Theorem \ref{new:thm}}. By means of the last Proposition and the special cases just treated, we have shown that the sections provided by one (hence both) of  the points $A,B$ are not identically torsion on the elliptic scheme over the conic $S:\mu^2=p+r\lambda^2$. Choose then the section $A$ for definiteness, By a well-known theorem of Silverman-Tate (see \cite{Z.UI}, Prop. 3.2 p. 69 or Appendix C by Masser)  the algebraic points $x$ of the conic $S$ such that $A_x$ becomes torsion (on the elliptic curve $\mathcal E_x$ corresponding to $x$) have bounded height. Therefore if we let $x$ vary though the rational points $S(\Q)$ of $S$, for all but finitely many ones $A_x$ will be non-torsion on $\mathcal E_x$. Therefore by taking multiples $nA_x, n\in\Z$ we shall obtain infinitely many rational points on $\mathcal E_x$ for such points $x$. Since $S(\Q)$ is non-empty by assumption, we have that $S$ is birationally equivalent to $\mathbb P_1$ over $\Q$, hence $S(\Q)$ is Zariski dense in $S$, whence the rational points are Zariski dense in the cubic family, concluding the argument.

\medskip

\subsection{An alternative construction}  
We present an alternative construction of a   section of the above cubic fibration.  It also provides rational points. 

 The equations  \eqref{eq:C_k} and \eqref{eq:C'k'}  and  give the following  linear relation between $x$ and $y$ with coefficients in $\Q[k,k']$
\begin{equation}
    \label{eq:phi}
    2p\varphi_1 x + 2 rk y = \varphi_0,
\end{equation}
where
\begin{equation}\label{eq:phi_1,2,3}
   \varphi_1 = k'-kc, 
   \qquad \varphi_0 = k'\delta- k\delta'.
\end{equation}
Then \eqref{eq:phi} allows to express $y$ as
   $ 2rky = \varphi_0-2p\varphi_1 x$.
Multiplying the equation $px^2+ry^2=z^2$ by $4rk^2$ and  substituting for $2rky$ from the last equation one gets
\begin{equation*}
4prk^2x^2+(2p\varphi_1x-\varphi_0)^2=4rk^2z^2;
\end{equation*}
using now the equation $2kz=2px-\delta$  one obtains 
\begin{equation}
    \label{eq:5}
    4prk^2x^2+(2p\varphi_1x-\varphi_0)^2=r(2px-\delta)^2\end{equation}
Hence  we get the following quadratic equation for $x$:
\begin{equation}
    \label{eq:ABC}
    A x^2 + B x + C = 0,
\end{equation}
where
\begin{equation}\label{eq:def-ABC}
\left\{ \begin{matrix}   A &=& 4prk^2+4p^2\varphi_1^2-4p^2r\\
    B &=&  4pr\delta-4p\varphi_0\varphi_1 \\
    C &=& \varphi_0^2-r\delta^2.
    \end{matrix}\right.
\end{equation}
We want to produce a rational solution to \eqref{eq:ABC}. One way to do this is to consider the case when $A = 0$, thus obtaining a linear equation in $x$ that has obviously a rational solution. We note that the condition $A=0$ amounts to
\begin{equation}
    p(k/p)^2+r(\varphi_1/r)^2=1.
\end{equation}
This represents a point on the conic $pX^2+rY^2=1$, which by assumption has a rational point, and thus can be parametrized rationally, say in  terms of a parameter $t$, as $pf(t)^2+rg(t)^2=1$, where $f,g$ are certain (quadratic) non-constant rational functions in $\Q(t)$. 
Putting $k=pf(t)$, $\varphi_1=rg(t)$,  we can express both $k$ and $k'$ as rational functions of $t$; in turn,   we can express  $\varphi_0$ as well rationally in terms of $t$.  Note that we can also use the former coordinates on the conic $S$, which would yield $rk\lambda=p\varphi_1$, $k=p/\mu$, whence   $k'=\varphi_1+kc=(r\lambda+cp)/\mu$.

The above formula then  yields an expression of $x=-C/B$ in terms of $t$ (or of $\lambda,\mu$), provided that $B$ is not identically zero as a function of $t$,  and therefore allows to express all the variables $x,y$ and $z$ as rational functions of $t$, using the above linear  equations  for $y$ and   for $z$. 

In order  to verify that $B$ is not identically  zero with this choice, we may argue  directly, or else note that in view of equation \eqref{eq:ABC} $B=0$ also implies $C=0$, since we are in the case $A=0$.   Hence, we have to exclude that $B\equiv C\equiv 0$ under  the present choice of $k,\varphi_1$. Note that $C=0$  amounts to  $\varphi_0^2=r\delta^2$  and under this condition $B=0$ implies 
$(r\delta)^2=r\delta^2\varphi_1^2$. But $\delta=p-k^2\neq 0$, whence $\varphi_1=r$ would be constant on the conic, which is not the case.
 
  This completes the verification that $x,y$, and so $z$, are well-defined in terms of $t$. 
Geometrically, this means that we have  defined a rational curve on the surface $\calS_3$, which can be viewed as a  section of our elliptic fibration. If we verify that this is not identically torsion, this section would  give rise to a Zariski-dense set of rational points, as in our former argument. We have not checked this property, since it appears a bit complicated, and since we have already done this for the former sections; however in principle one can argue as above, exploring for instance when the   section becomes of order $2$, or of any chosen order, and trying to prove that this happens outside the bad reduction. 

\section{The Geometry of the surface \texorpdfstring{$\calS_3$}{S3}}\label{sec:S3geo}
The goal of this section is to study the geometry of the surface 
\[
    \calS_3 = \calS(O,P,P') = \{ Q \in \A^2(\Q) : d(O,Q)\in \Q,\ d(P,Q) \in \Q \text{ and } d(P',Q) \in \Q \},
\]
i.e. the set of points $Q$ with rational distances from $O,P$ and $P'$. For this analysis we restrict to the case in which $P,P' \in \Q^2$ so that any point with rational distance from $O,P$ and $P'$ is automatically defined over $\Q$ (provided that the points are not aligned, which we have already discussed in a previous section).

A rational point $Q = (x,y) \in \calS_3(\Q)$ corresponds to a rational solution to the system
    \begin{align}
        &x^2 + y^2 = z^2,  \label{eq:1uv} \\ 
        &(x - a)^2 + (y - b)^2 = u^2, \label{eq:2uv} \\
        &(x - a')^2 + (y - b')^2 = v^2, \label{eq:3uv}
    \end{align}
where $z = d(O,Q)$, $u = d(P,Q)$ and $v = d(P',Q)$. In this presentation  the set $\calS_3$ lies in the intersection of $3$ quadrics in $\A^5$,  defined by the three equations above. We shall show that (a smooth projective completion of) this variety is a $K3$ elliptic surface with a dense set of rational points over $\Q$.

In the sequel we shall also prove that the surface $\calS_3$ is endowed with (several) elliptic fibrations parametrized by a line. The presence of a second fibration is related to the density of rational points even in the sense of the Euclidean topology, as first noticed by Swinnerton-Dyer in \cite{Swinnerton}. 

\subsection{The surface \texorpdfstring{$\calS_3$}{S3} as a ramified cover of \texorpdfstring{$\PP_2$}{P2}}\label{SS.geom}
Let us show how to re-obtain the surface $\calS_3$ as an abstract ramified covering of the $(x,y)$-plane. In the first part of this paragraph we argue birationally, describing the surface     $\calS_3$ only up to birational transformations. In the last part of the paragraph, we shall work out a smooth complete model.

We first note that the three equations above correspond each to a quadratic cover of the plane ramified over the union of two lines: for instance, \eqref{eq:1uv} amounts to taking a cover ramified over the singular conic of equation $x^2+y^2=0$, which consists of  a pair of complex conjugate lines.

The two other pairs of lines, defined by the vanishing of the left-hand side in equations \eqref{eq:2uv}, \eqref{eq:3uv}, have the same slopes  of the first pair; hence the six lines intersect in triples at the two points at infinity $(1:\pm i:0)$ (under the embedding $\A^2\hookrightarrow \PP_2$ given by $(x,y)\mapsto (x:y:1)$).
\medskip

As we said, the surface $\calS_3$ is obtained as a degree-eight covering of the plane $\PP_2$, actually an abelian covering of type $(2,2,2)$. This covering can then be decomposed as
$$
\calS_3\to X_2\to X_1\to \PP_2
$$
where each arrow denotes a degree-two covering. (Note that we can choose the surface $X_1$ to be birationally equivalent to the surface $\calS_1(O)$ and $X_2$ to be birationally equivalent to the surface $\calS_2(O,P)$.) 

For instance, for the first covering $X_1\to\PP_2$ we can take  the one ramified over the pair of lines $x^2+y^2=0$; the second covering is algebraically described by adjoining to the function field of $X_1$ the square root of the function $(x-a)^2+(y-b)^2$; it ramifies  over the zero set of this function, i.e. the pre-image in $X_1$ of the mentioned pair of lines intersecting on $P$;  the last one ramifies over the pre-image in $X_2$ of the pair of lines intersecting on $P'$. 

  Let us consider the pencil $\Lambda_0$ of lines in $\PP_2$ passing through the origin $0$ of $\A^2\subset \PP_2$. Each such line has a reducible pre-image in $X_1$, consisting of two components. Each such component $l$ is a smooth  curve of genus zero above which the cover $\calS_3\to X_1$ induces a cover of type $(2,2)$. Since the cover $\calS_3\to X_1$  ramifies above the four curves in $X_1$ defined by the vanishing of the functions $(x - a)^2 + (y - b)^2=0$ and $(x - a')^2 + (y - b')^2=0$, and $l$ intersects each such curve at a single point, the ramification on $l$ occurs above four points. Hence, by Riemann-Hurwitz, 
   the  pre-image in $\calS_3$ of each such curve $l$ is then a genus-one curve. In this way, we obtain a family of genus-one curves on  $\calS_3$.  Changing the role of the three equations \eqref{eq:1uv}, \eqref{eq:2uv}, \eqref{eq:3uv} one can consider the pencils $\Lambda_P$ (resp. $\Lambda_{P'}$) of lines through $P$ (resp. $P'$), obtaining a total of three families of genus-one curves.
   
Let us consider the projection $\calS_3 \dashrightarrow \Lambda_0$. It is defined by projecting a point  $s=(x,y,z,u,v)\in \calS_3$ to the corresponding point $(x,y)$ of the plane and then associating the line joining the origin to the point $(x,y)$ (which can be identified with the projective point $(0:0)$. As mentioned, the pre-image of any line in $\Lambda_0$  is reducible, unless such line is a component of the singular conic $x^2+y^2=0$; in that case, the pre-image is a double curve (a non-reduced divisor). Hence there cannot exist sections $\Lambda_0\to \calS_3$ of the above defined projection $\calS_3 \dashrightarrow \Lambda_0$. However, considering the Stein factorization
\begin{equation}\label{E.fibrazione1}
\calS_3 \dashrightarrow \PP_1 \to \Lambda_0
\end{equation}
where $\PP_1\to \Lambda_0$ is the quadratic cover ramified over the two points $(1:\pm i)$, the first projection 
 \begin{equation}\label{E.fibrazione2}
\pi_0: \calS_3\dashrightarrow \PP_1
\end{equation} 
which is a fibration (i.e. its fibers are generically irreducible) 
does admit sections, as we shall show below.

\smallskip

Let us now describe a smooth projective model of the surface $\calS_3$ which regularizes the projection $\calS\dashrightarrow \Lambda_0$. By abuse of notation, we shall denote again by $\calS_3$ such a model. To construct this model, we shall proceed as in section \ref{S.Intermezzo}. We start from blowing-up the three points $O,P,P'$ on $\PP_2$, together with the two points at infinity $(1:\pm i:0)$. We obtain a smooth surface $\calS$. The surface $\calS_3$ is then defined as the $(2,2,2)$ covering of $\calS$ ramified over the strict transform of the singular conics defined by the quadratic forms appearing on the left-hand side of \eqref{eq:1uv}, \eqref{eq:2uv}, \eqref{eq:3uv}.  More precisely, $\calS_3$ will be the normalization  of the fiber product of the three surfaces each of which obtained as the quadratic cover of $\calS$ ramified over each of the mentioned degenerate conics.

A simple (but lengthy) application of the Riemann-Hurwitz formula for surfaces shows that the canonical class of $\calS_3$ vanishes, so the surface is a $K3$ surface (or an abelian surface, but this possibility is ruled out by the fact that it contains  non-constant families of elliptic curves).
Alternatively, one can consider a generic line on $\PP_2$: its pre-image on $\calS_3$ is a curve $\mathcal C$ covering the line eight-to-one; more precisely, $\mathcal{C}\to\PP_1$ is an abelian cover of type $(2,2,2)$ ramified over eight points; then its Euler characteristic is $8$; since $\mathcal{C}^2=8$, by adjunction it follows that the product $\mathcal{C}\cdot K_{\calS_3}$ vanishes, where $K_{\calS_3}$ denotes a canonical divisor of $\calS_3$. Hence the canonical class is numerically trivial.

\subsection{Alternative models for \texorpdfstring{$\calS_3$}{S3}} Let us consider once again the diagram \eqref{E.fibrazione1} leading to the projection $\pi_O:\calS_3\to \PP_1$ of equation \eqref{E.fibrazione2}. Changing $O$ with the two other points $P,P'$ leads to a rational map 
$$
\calS_3 \dashrightarrow \PP_1\times\PP_1\times\PP_1.
$$
Regularizing this map, and abusing again of the notation by using again the symbol $\calS_3$ for the new domain, one obtains a regular map  of $\calS_3$ into $\PP_1\times\PP_1\times\PP_1$, actually an embedding. We now show that the image has multi-degree $(2,2,2)$. This amounts to saying that each of the three projections $\calS_3\to (\PP_1)^2$, obtained by omitting one factor in the projection to $(\PP_1)^3$, has degree two. Let us prove this fact. The projection $\calS_3\to \PP_1\times\PP_1$ can be composed on the left to obtain a projection 
 \begin{equation}\label{E.proiezione8}
 \calS_3\to \Lambda_P \times \Lambda_{P'}
 \end{equation}
  (see again diagram \eqref{E.fibrazione1}), whose degree is four times the degree of the map $\calS_3\to (\PP_1)^2$ we are interested in. Now, the fiber of a generic point in $\Lambda_P\times\Lambda_{P'}$ is a single point in $\PP_2$ (since a point in $\Lambda_P$ (resp. in $\Lambda_{P'}$) represents a line in $\PP_2$ and two lines intersect at one point) and so consists in eight points in $\calS_3$. Hence the projection \eqref{E.proiezione8} has degree eight and so the projection $\calS_3\to (\PP_1)^2$ has degree two as wanted.
 
 Now, smooth surfaces in $\PP_1\times\PP_1\times\PP_1$ of multi-degree $(2,2,2)$ are known to be $K$3 surfaces. Actually, the surface studied in \cite{CZ} by the first and the third author in the frame of the  Hilbert Property belongs to this family of $K$3 surfaces.
 
\smallskip

\subsection{The Kummer model for the surface  \texorpdfstring{$\calS_3$}{S3}}\label{subsec:kummer} In this section we shall realize the surface $\calS_3$ as a quartic in $\PP_3$. First note that working over the the complex field (and actually even over the reals) we can choose coordinates so that the second point $P$ is the point $(1,0)$; of course, for our Diophantine problem we must then consider twists of the surface we are constructing.  The three equations defining $\calS_3$ become then
 
\begin{align}
        &x^2 + y^2 = z^2   \nonumber \\ 
        &(x - 1)^2 + y^2 = u^2   \label{E.S3Speciale} \\
        & (x-a)^2 + (y - b)^2 = v^2 \nonumber, 
    \end{align}
    where $a\in \R$ and $b\neq 0$. 
Using the first equation to eliminate $x^2,y^2$ from the second and third equations, we obtain
\[
2x=z^2-u^2+1,\qquad 2ax+2by=a^2+b^2+z^2-v^2,
\]
hence 
\begin{align*}
        &x   = \dfrac{z^2-u^2+1}{2}   \\ 
        &y  = \dfrac{a^2+b^2+(1-a)z^2-v^2+au^2-a}{2b} \\
    \end{align*}
    which together with the first equation \eqref{E.S3Speciale} give rise to a single quartic equation in $z,u,v$. For instance 
when $a=0, b=1$ we obtain the following equation, that was already present in \cite[Section D19]{Guy}:
\begin{equation}\label{E.singola}
2z^4+v^4+u^4+2=2z^2(v^2+u^2)+2u^2+2v^2.
\end{equation}
The surface defined by this equation, and by those obtained for arbitrary choices of $a,b$ with $ b\neq 0$, is a singular quartic surface, having sixteen simple nodes . 

Recall that a Kummer surface is a smooth projective minimal model of a quotient of the form $A/\{\pm 1\}$, where $A$ is an abelian surface and $\{\pm 1\}$ denotes the group generated by the map $A\ni p\mapsto -p\in A$; sometimes one calls a Kummer surface also such a quotient, which has   sixteen isolated singularities, corresponding to the $2$-torsion of the abelian surface $A$. The corresponding smooth model has sixteen $-2$-curves, projecting to such singular points.  

We shall denote by $\mathrm{Kum}(A)$ the Kummer surface associated to the abelian surface $A$.

\begin{theorem}\label{T.Kummer}
The surface $\calS_3$ is a Kummer surface.
\end{theorem}

It is a classical fact, see for example \cite{Nik}, that every complex quartic surface on $\PP_3$ with sixteen nodes, and no other singularities, is a (singular model of a) Kummer surface; such a surface is actually isomorphic to a surface of the form $A/\{\pm 1\}$, where the sixteen nodes correspond to the sixteen points of $2$-torsion on $A$, which are fixed by the involution on $A$.

 This would suffice to  prove the result. However, we shall give an explicit proof of Theorem \ref{T.Kummer},  constructing  the abelian surface $A$ which covers $\calS_3$ and describing the cover $A\to \calS_3$.  

\begin{proof}
As we said, over the complex numbers we can suppose the equations of $\calS_3$ are those given in \eqref{E.S3Speciale}. After setting
$$
x+iy=\xi,\qquad x-iy=\eta
$$
the three equations \eqref{E.S3Speciale} take the form
\begin{align} 
        &\xi\eta = z^2    \\ 
        &(\xi - 1)(\eta-1) = u^2   \\
        & (\xi-a-bi)(\eta-a+bi) = v^2, 
    \end{align}
Consider now the two elliptic curves $E_1$ and $E_2$ given in affine equations as
\begin{equation}\label{E.elliptic-curves}
E_1:\, y_1^2=\xi(\xi-1)(\xi-a-bi),\qquad  E_2:\, y_2^2= \eta(\eta-1)(\eta-a+bi).
\end{equation}
We note at once that $E_1,E_2$ are not defined over $\Q$ and they contain only finitely many $\Q$-rational points (in fact only the points of order $2$): in particular the rational points on the associated Kummer surface $\mathrm{Kum}(E_1\times E_2)$, do not come from the rational points on $E_1\times E_2$.

The function field of the surface $\calS_3/\Q$ contains the square root of the product 
$$\xi(\xi-1)(\xi-a-bi)\eta(\eta-1)(\eta-a+bi)=(x^2+y^2)((x-1)^2+y^2)((x-a)^2+(y-b)^2);$$ this square root  together with the functions $x$ and $y$ generate  the function field (even over $\Q$) of the Kummer surface $\mathrm{Kum}(E_1\times E_2)$. 
In particular, after putting $\alpha=a+bi\in\C\setminus \R$, we obtain 
$$
\C(\mathrm{Kum}(E_1\times E_2))=\C(\xi,\eta)(\sqrt{\xi(\xi-1)(\xi-\alpha)\eta(\eta-1)(\eta-\bar{\alpha})})
$$
and the inclusion of functions fields
\begin{align*}
\C(\xi,\eta)&\subset \C(\xi,\eta)(\sqrt{\xi(\xi-1)(\xi-\alpha)\eta(\eta-1)(\eta-\bar{\alpha})}) \\
&\subset\C(\xi,\eta)(\sqrt{\xi(\xi-1)(\xi-\alpha)},\sqrt{\eta(\eta-1)(\eta-\bar{\alpha})})=\C(E_1\times E_2).
\end{align*}
Geometrically, this corresponds to the dominant maps 
\begin{equation}\label{E.KUm1}
\PP_1 \times \PP_1 \leftarrow \mathrm{Kum}(E_1\times E_2)\leftarrow E_1\times E_2. 
\end{equation}
We also have a rational dominant map of degree $4$, actually an abelian cover of type $(2,2)$,
$$
\calS_3\to \mathrm{Kum}(E_1\times E_2).
$$
This map corresponds to the inclusion
\begin{equation*}
\C(\xi,\eta)(\sqrt{\xi(\xi-1)(\xi-\alpha)\eta(\eta-1)(\eta-\bar{\alpha})})  \subset \C(\xi,\eta)(\sqrt{\xi\eta},\sqrt{(\xi-1)(\eta-1)},\sqrt{(\xi-\alpha)(\eta-\bar{\alpha})})
\end{equation*}
Consider now the degree-$2$ map $E_1\times E_2\to \mathrm{Kum}(E_1\times E_2)$; performing the fiber product of the two maps pointing to $ \mathrm{Kum}(E_1\times E_2)$ we obtain a commutative diagram, where $A\to \calS_3$ is a degree-$2$ map and the arrow $\phi: A\to E_1\times E_2$ is a degree-$4$ Galois covering of type $(2,2)$. 
\begin{equation}\label{diagram}
   \xymatrix{A \ar[r]^<<<<<<\phi \ar[d] & E_1 \times E_2 \ar[d] \\
   \calS_3 \ar[r] \ar[rd] & \mathrm{Kum}(E_1 \times E_2) \ar[d] \\
   & \PP_1 \times \PP_1} 
\end{equation}

We now prove that $\phi$ is unramified, thus proving that $A$ is an abelian surface. 

Actually, this cover corresponds to the function field extension
$$
\C( E_1\times E_2)\subset \C(\xi,\eta)(\sqrt{\xi(\xi-1)(\xi-\alpha)},\sqrt{\eta(\eta-1)(\eta-\bar{\alpha})}, \sqrt{\xi\eta }\sqrt{(\xi-1)(\eta-1)}). 
$$
Now, if we look at the possible ramification of this extension, this can occur only over the curves on $E_1\times E_2$ lying over the lines $\xi=0, \eta=0, \xi-1=0, \eta-1=0$ of $\PP_1\times \PP_1$; moreover the ramification index can be at most $2$. However, these lines are already ramified with index $2$ on the cover $\calS_3\to \PP_1\times \PP_1$; by Abyankhar's lemma, the cover $A\to E_1\times E_2$ turns out to be unramified.
\end{proof}
\bigskip 

\noindent{\tt Weak approximation.} A natural question, also  related   to the distribution of rational points in $\calS_3$,  is whether weak-approximation holds. This has been an extensively studied topic both for elliptic K3 surfaces and for Kummer surfaces. For example, there are known instances of elliptic K3 surfaces for which there is a `transcendental obstruction' to Weak Approximation, in the sense of the paper \cite{HVV} . More generally Colliot-Th\'el\`ene, Skorobogatov and Swinnerton-Dyer in \cite{CTSSD}, building on previous work of Swinnerton-Dyer \cite{SD}, studied the Brauer-Manin obstruction to weak approximation for elliptic surfaces, and Skorobogatov and Swinnerton-Dyer analyzed, in particular, the case of a family of Kummer surfaces, whose associated abelian variety is the product of two elliptic curves. They proved in \cite[Theorem 1]{SSD} that, under certain hypotheses (for example that the Jacobians of the two elliptic curves have all their 2-division points defined over the ground field, plus additional more technical conditions),  the equations defining these families of Kummer surfaces are soluble in the base field $k$.

\smallskip
 
\noindent{\tt Field moduli and field of definition}. As we remarked, the condition $(ii)$ about the points $P_1,P_2,P_3$ in Theorem \ref{T.main} is equivalent to saying that the relevant (Kummer $K3$) surface $\calS_3=\calS_{P_1,P_2,P_3}$ can be defined over $\Q$. The following example shows that such surfaces can have the rational field $\Q$ as a field of moduli without being definable over $\Q$. Consider the triangle $0,(1,1\sqrt{2}),(1,1-\sqrt{2})$: it is defined over $\Q(\sqrt{2})$ and is invariant by Galois conjugation. Hence the corresponding surface is defined over $\Q(\sqrt{2})$ and is isomorphic to its Galois conjugate. However, the quadratic form appearing in $(ii)$ of Theorem \ref{T.main} is not defined over $\Q$, so the surface is not definable over $\Q$.

\smallskip

\noindent{\tt The Hilbert property}. The rational points in the surfaces of type $\calS_3$ turn out also to satisfy another density condition, stronger then Zariski-density and not comparable with density in the real topology. Namely they satisfy the Hilbert Property; equivalently, they form a {\it non-thin set}. The Hilbert property for such surfaces   can be derived, after proving the Zariski-density of rational points and non-torsion of the sections of the elliptic fibrations considered in this paper, by the methods introduced in \cite{CZ} and developed by J. Demeio in \cite{Demeio}, \cite{Demeio2}.

\section{Distance from four or more points}\label{sec:4pts}
 Let us now consider the problem of the rationality of the distances from four (or more) points. This corresponds to adding a quadratic equation of the form
$$
(x-\alpha)^2+(y-\beta)^2=w^2
$$
to the system of the three equations  \eqref{eq:1}, \eqref{eq:2}, \eqref{eq:3}. Geometrically, adding the extra equation  corresponds to operating a quadratic cover $\calS_4\to\calS_3$  of the $K3$-surface  $\calS_3$ ramified over the pull-back of the curve defined by the vanishing of the left-hand side of the above equation. This equation represents again a pair of lines on the plane, and a two-component curve on $\calS$ of self-intersection $2$.
This is a big (but non-ample) divisor on $\calS$.  Its pull-back to $\calS_3$ is again a big divisor on $\calS_3$, say $D$. Then any cover of the $K3$-surface $\calS_3$ ramified over $D$ is a surface of general type.  

According to the celebrated Bombieri-Lang conjecture, its rational points (over any fixed number field) should be degenerate. In particular we have the following fact:

\begin{proposition}\label{prop:S4_BL}
   The Bombieri-Lang conjecture implies that the $k$-rational points on the surface $\calS_4$ are not Zariski dense, for any number field $k$.
\end{proposition}

However, in our situation, we lack techniques to study the degeneracy of rational points. In particular there is no known case of the conjecture for non-singular simply connected surfaces\footnote{see \cite{GRTW} for a recent example over function fields.} and $\calS_4$ is indeed simply connected.

\begin{proposition}\label{prop:SimpConn}
    The surface $\calS_4$ is simply connected.
\end{proposition}
\begin{proof}
    Using the same arguments as in Section \ref{sec:S3geo} we can describe the surface $\calS_4$ as a ramified cover of the blow up of $\PP_1 \times \PP_1$ over (the strict transform of) 5 horizontal lines and 5 vertical lines, of type $(2,2,2,2)$. To conclude that the surface is simply-connected one can use the description of the fundamental group of a bidouble cover given in \cite{Cat} (see for example Proposition 2.7 or more generally Proposition 1.8).
\end{proof}

An alternative approach is via the so-called \emph{product-quotient surfaces} studied by Catanese et al., where the surface $\calS_4$ has the same fundamental group as the quotient of a product of two curves of genus 5 (which one can construct explicitly from the map $\calS_4 \to \PP_1 \times \PP_1$). The latter quotient has trivial fundamental group using the main Theorem of Armstrong in \cite{Arm}. Indeed, using the description of $\calS_4$ as a $(2,2,2,2)$ Galois cover of the blow up $\PP_1 \times \PP_1$, with group $(\Z/2\Z)^4$, extending the similar construction for $\calS_3$ given in Section \ref{sec:S3geo}, one obtains two genus $5$ curves $C_1,C_2$, both covers of $\PP^1$ of type $(2,2,2,2)$ (these are constructed using the 5 vertical, and 5 horizontal lines in the branch divisor of the cover $\calS_4$). Their product, $\calC_1 \times \calC_2$, dominates $\calS_4$, and satisfy the hypotheses of \cite{Arm}, thus obtaining that $\calS_4$, is simply connected.

\medskip

We end this section with an easy observation: a proof of the fact that $\calS_4(\Q)$ is not Zariski dense, would answer in particular  the famous Erd\H{o}s-Ulam problem in the negative, i.e. there would be no Zariski-dense rational distance set in the plane. Recall that a rational distance set is a set of points whose mutual distances are all rational. Given such a set $X\subset\R^2$,   after choosing four points  $P_1,\ldots,P_4$ in general position in $X$, supposing that $X$ does contain at least four points in general position,  the set $X\setminus\{P_1,\ldots,P_4\}$ can be lifted to a set of rational points on a surface $\calS_4$. The conjectural degeneracy of the set  $\calS_4(\Q)$ immediately implies the degeneracy (in the Zariski topology) of the set $X$. In turn, this last fact easily implies that all but finitely many points of $X$ must be contained in a line or a circle.

\medskip
\section{Integral Distances from two points}

As mentioned above, in \cite{Z} it was shown that the {\it integral} (=lattice) points of the plane having integral distance from two lattice points $O,P$ are {\it usually} infinite in number, but are never Zariski dense (and in fact always contained in a finite union of hyperbolae and lines). This last fact may be read as  a case of Runge's theorem for surfaces.   For larger rings of integers the situation is different (as it happens for Runge's theorem for curves):  in fact, as we now show,   over certain other  rings of algebraic integers, the integral points in question form a Zariski-dense set.  

In this section we study solutions to the system given by equations \eqref{eq:1_2pt} and \eqref{eq:2_2pt} in rings of integers of certain  {\it quadratic} number fields.  
We shall work out a couple of explicit cases, both real and non-real. 
In that last case then the meaning of `distance' loses its usual significance.

 \subsection{Integral points with coordinates in a real quadratic number field}
 Let $\delta>0$ be any positive integer, not a perfect square. Then the quadratic ring $\Z[\sqrt{\delta}]$ admits infinitely many units. This fact is crucial in the proof of the following result:
 
 \begin{theorem}\label{T.quadratico-reale}
 Let $O,P\in\Z^2$ be integer points of the plane and let $\delta>0$ be a  non-square positive integer. The set of points with coordinates in the ring $\Z[\sqrt{\delta}]$ whose distances from $O$ and from $P$ belong  to the ring $\Z[\sqrt{\delta}]$ is Zariski-dense in the plane.
 \end{theorem}

\begin{proof}

Choosing coordinates, as before, such that $O=(0,0)$ and $P=(a,b)$, consider the line $l$ of equation $x=a$. For every point $Q=(a,t)\in l$ with $t\in \Z[\sqrt{\delta}]$, its distance from  $P$ equals $|t-b|$, hence lies in the quadratic ring  $\Z[\sqrt{\delta}]$. The condition that $d(O,Q)$  also belongs to that ring amounts to saying that 
$$
a^2+(b-t)^2=s^2,
$$
 for some $s\in \Z[\sqrt{\delta}]$. This is the equation of a conic admitting the integral points $(t,s)=(b,a)$. Due to the presence of infinitely many units in the ring $\Z[\sqrt{\delta}]$, that conic admits infinitely many points with coordinates in that ring. This is a well-known fact, however we give a short proof for simplicity.  
 Let us write the above equation as
 $$
 (s-b+t)(s+b-t)=a^2.
 $$
 For every unit $\omega\in(\Z[\sqrt{\delta}])^*$,   solving for $s,t$ the equations 
$$
s-b+t=a\omega,\qquad s+b-t=b-a\omega^{-1}
$$
provides a rational solution. This solution is integral whenever $\omega\equiv 1\pmod 2$, which happens for an infinite subgroup of units.

Hence the line $l$ contains an infinite set $X\subset l$ of $\Z[\sqrt{\delta}]$-integral points whose distances  from $O,P$ belong to $\Z[\sqrt{\delta}]$. The same of course holds for the horizontal line  $y=0$. To produce a Zariski-dense set we argue as in section \ref{sec:2pt_conics}. 

Consider the family of confocal conics with foci $O$ and $P$. Now, two such conics pass through every point   $x\in X$. Let us choose e.g. the hyperbola $\mathcal{C}_x$  passing through $x$; as in the previous argument, such  hyperbola $\mathcal{C}_x$ contains infinitely many integral points (with respect to the ring $\Z[\sqrt{\delta}]$). For each such integral point $y$ we have  that the three quantities  $d(O,y)^2,d(P,y)^2$ and $d(O,y)-d(P,y)$ are all integral. Then $d(O,y)+d(P,y)$ is rational; since it is an algebraic integer (in view of the fact that $d(O,y)^2,d(P,y)^2$ belong to $\Z[\sqrt{\delta}]$), it is also an element of $\Z[\sqrt{\delta}]$. It follows that the integral points of $\mathcal{C}_y$ are still solutions to our problem. Since the union of the $\mathcal{C}_x$, $x\in X$, is Zariski-dense, the theorem is proved.
\end{proof}

It might be interesting to look at the points with coordinates in the ring $\Z[\sqrt{\delta]}$ having $\Q$-{\it rational}   distance from $O$ and $P$. However, in this case the distances turn out to be rational integers (in $\Z$), and by a Runge-type argument (as in the mentioned case of integral points with integral distances) the degeneracy follows.

A more compelling problem consists in studying  the distribution of the points $Q$ having coordinates in the number field $\Q(\sqrt{\delta})$ and having $\Q$-rational distances from both $O$ and $P$. 
This amounts to intersecting circles centered in $O$ and $P$ with rational radii and asking that the intersection points, which generically are quadratic over $\Q$, belong to a fixed number ring $\Q(\sqrt{\delta})$. Of course this set of points contains the set $\calS_2(\Q)$ considered in the previous section, and consisting in rational points.  We shall then concentrate on those points which are irrational but have rational distances from the given points.

Writing 
$$
x=x_1+\sqrt{\delta}x_2,\qquad y=y_1+\sqrt{\delta}y_2,
$$
for rational numbers $x_1,x_2,y_1,y_2$, the condition that $d(O,Q)^2$, where $Q=(x,y) $,  is rational leads to the relation
$$
x_1x_2+y_1y_2=0.
$$
We shall neglect the already mentioned case $x_2=y_2=0$ (i.e. $Q$ rational), as well as the case $x_iy_j=0$ for any choice of $i,j\in \{1,2\}$. The rationality of the squared distance $d(P,Q)^2$ gives an analogous relation. The two relations together lead to an expression of $y$ in terms of $x$, namely:
$$
y_1=\frac{b}{a}x_1,\qquad y_2=-\frac{a}{b} x_2.
$$
Then the rationality of the distance $d(O,Q)$ gives
$$
x_1^2 \left(1+\frac{b^2}{a^2}\right)+\delta x_2^2\left(1+\frac{a^2}{b^2}\right)=z^2,
$$
which can be re-written, after renaming the variables, as
$$
u^2+\delta v^2=(a^2+b^2)w^2.
$$
It is clear that for certain choices of $\delta,a,b$, the above equation has no non-trivial rational solution. The rationality of the distance $d(P,Q)$ provides another (non-homogeneous) quadratic equation. The two equations together define again a Del Pezzo surface of degree four, as in the previous section. Given the points $P=(a,b)$, for some choice of $\delta$   the set of points defined on $\Q(\sqrt{\delta})$, but not over $\Q$, and having rational distance from both $O$ and $P$ will be dense.

\medskip

\subsection{Points with coordinates in \texorpdfstring{$\Z[i]$}{Z[i]}} 
Let us  now consider an imaginary quadratic  ring. 
The simplest such case occurs when  both the coordinates of the points $Q$ and the `distances' $d(O,Q), d(P,Q)$ are Gaussian integers, i.e. elements of $\Z[i]$. In this setting, we keep the {\it algebraic } definition of  distance, namely the distance between  two points $A,B\in (\Z[i])^2$ is  a square-root of the 
sum of the squares of the entries of the vector $\vec{AB}$; this is defined only up to sign, but that is immaterial for the field of definition. Basically, we are interested on the points of the surface $\calS_2$ defined over the ring $\Z[i]$. 
We will show that such set  is Zariski dense:

\begin{proposition}\label{prop:Zi}
   Given $P=(a,b) \in \Z^2$, $P\neq O$,  the set $ \calS_2(\Z[i])$ is Zariski dense, where as usual
   \[
    \calS_2(\Z[i]) = \{ Q \in \A^2(\Z[i]) : d(O,Q)\in \Z[i], d(P,Q) \in \Z[i] \}.
   \] 
\end{proposition}

\begin{proof}
    
We start by observing that, similarly to the case of  rational distances, the set $\calS_2(\Z[i])$ corresponds to the system given by equations \eqref{eq:1_2pt} and \eqref{eq:2_2pt}. Since now $-1$ is a square in our ring, the system is equivalent to the following: 
\begin{align}
    &w_1   w_2 = z^2,  \label{eq:i1} \\ 
    &(w_1 - \alpha) (w_2 - \bar\alpha)= (z - k)^2 \label{eq:i2},
\end{align} 
where we denote with $\bar\alpha $ the number $a-ib$ (which is the complex conjugate of $\alpha$ if $a,b\in\R$, as we are assuming). By unique factorization in the ring $\Z[i]$, we can rewrite the system as follows: we set $w_1 = d_1 r^2$ and $w_2 = d_1 u^2$, and similarly $w - \alpha = d_2 s^2$ and $ w_2 - \bar\alpha = d_2 v^2$ for some $d_1,d_2,r,s,u,v \in \Z[i]$ (where for instance $d_1$ is a $\gcd$ of $w_1,w_2$). Every time we can write $w_1,w_2,w_1-\alpha, w_2-\bar{\alpha}$ in such a way for suitable $(d_1,d_2,r,s,u,v)$ we obtain a solution of the system of equations \eqref{eq:i1}, \eqref{eq:i2}. Eliminating $w_1,w_2$ this amounts to the following system:
\begin{align}
    &d_1 r^2 - d_2 s^2 = \alpha, \label{eq:i3} \\ 
    &d_1 u^2 - d_2 v^2 = \bar\alpha \label{eq:i4}. 
\end{align}

We want to show that the above system has infinitely many solutions $(d_1,d_2,r,s,u,v)\in(\Z[i])^6$,  and that these solutions give rise to a Zariski-dense set of points $(w_1,w_2)=(d_1r^2,d_1u^2)\in\A^2$; of course, the density of the pairs $(w_1,w_2)$ is equivalent to the density of the pairs $(x,y)$, which are obtained from $(w_1,w_2)$ via a fixed linear transformation. 

This amounts to show that for given $\alpha\in\Z[i], \alpha\neq 0$  there exists choices of non-zero $d_1,d_2$ such that each of the equations \eqref{eq:i3}, \eqref{eq:i4} has infinitely many solutions. Indeed, in view of the fact that the two equations have separated variables, if each of two equations \eqref{eq:i3}, \eqref{eq:i4} has infinitely many solutions  the pairs $(r,u)$ produced  from these solutions will be Zariski-dense on the plane; consequently also the pairs $(w_1,w_2)=d_1(r^2,u^2)$ will be Zariski-dense.  

The form of the two equations above is amenable to be studied using the theory of Pell's equations.

One can see from the above system that there always exist a choice of $d_1,d_2$ such that the two equations \eqref{eq:i3} and \eqref{eq:i4} have both a solution: it is enough to evaluate the left hand sides at a specific choice of $r,s,u$ and $v$ and solve for $d_1$ and $d_2$ (for example choosing $r,s,u,v$ that form a matrix with determinant invertible in $\Z[i]$). Therefore we are in the situation in which we have two conics defined over $\Z[i]$ both with an integral point over $\Z[i]$. However, such conics might contain only finitely many points with coordinates in $\Z[i]$. 

To avoid this case,  we can use (a case of) a theorem of Alvanos, Bilu and Poulakis \cite{ABP} that, in this special case, states the following: for an affine curve $\calC$ of genus zero, defined over a number field $K$, to have infinitely many integral points is sufficient to have two points at infinity which are conjugate over $K$ and whose field of definition is not a CM-extension of $K$ (see also \cite{Silv} for the analogue statement over the integers, due to Gauss). Recall that an extension of number fields $L/K$ is called a \emph{CM-extension} if $K$ is totally real and $L$ is totally imaginary. In our case, there are two points at infinity of equations \eqref{eq:i3} and \eqref{eq:i4}, namely $[\pm\sqrt{D}:1:0]$, where $D = d_2/d_1$. Hence our system will have infinitely many solutions provided that $d_2/d_1$ is not a square in $\Q[i]$. (Note that this condition is relevant: in fact, Theorem 1.2 in \cite{ABP} shows that when $\sqrt{d_2/d_1} \in \Q(i)$ the integral points will be finite in number. Note also that, when $\sqrt{d_2/d_1} \notin \Q(i)$, since $\Q(i)$ is totally imaginary, the field of definition of the points at infinity is never a CM extension of $\Q(i)$.)

Finally, it it remains to prove the following 

\begin{lemma}
Given a non-zero number $\alpha\in\Z[i]$ there exist $d_1,d_2,r,s,u,v\in\Z[i]$ such that $d_1/d_2$ is not a square and the relations \eqref{eq:i3}, \eqref{eq:i4} hold.
\end{lemma}

\begin{proof}
Let us  separate the proof into two cases, according as  $\alpha\bar\alpha$ being a perfect square or not.

If $\alpha\bar\alpha$ is not a square, which is equivalent to saying that $\alpha/\bar\alpha$ is not a square in $\Q(i)$, then take $d_1=\alpha, d_2=\bar{\alpha}$, $r=1, s=o, u=0, v=i$.

If, on the contrary $\alpha\bar\alpha$ is a square, by unique factorization in $\Z[i]$ we can write
$$
\alpha=\rho \cdot \xi^2, \quad \tilde\alpha=\bar\rho \cdot {\bar\xi}^2,
$$
where $\rho$  is square-free. The decomposition is unique up to units. Also, the fact that $\alpha\bar\alpha$ is a perfect square forces $\rho\bar\rho$ to be a perfect square; then $\rho$, being square-free,  cannot be divisible by non-rational primes; it follows that $\rho$ is real (hence a rational integer) or purely imaginary. In the first case, choose $d_1,d_2\in\Z$ with $d_1-d_2=\rho$ and $d_1/d_2$ not a square, which is certainly possible. Then put $r=s=\xi$, $u=v=\bar{\xi}$. In the second case, i.e. when $\rho$ is purely imaginary, take two purely imaginary numbers $d_1,d_2$ with the above property and $r=s=\xi$, $u=v=i\bar{\xi}$. This completes the proof of the lemma.
\end{proof}
 
To conclude the proof of the proposition, note that given a point $P(a,b)\neq (0,0)\in\Z^2$, the complex number $\alpha=a+bi\in\Z[i]$ is non-zero. Then by the above lemma and the previous considerations the system \eqref{eq:i3}, \eqref{eq:i4} has infinitely many solutions in $\Z[i]$ giving rise, as explained, to a dense set of points $Q\in \Z[i]\times \Z[i]$ in the plane with $d(O,Q), d(P,Q)$ both in $\Z[i]$ as wanted. 
\end{proof}

As we mentioned in the introduction, the set of  points having integral distances from three or more fixed points, over any ring of $S$-integers in any number field, is expected to be degenerate. Indeed, it is (the projection to the plane of) the set of integral points on an affine K3 surface; since the canonical divisor of a K3 surface is zero, the sum of the divisor at infinity and the canonical divisor turns out to be ample and Vojta's conjecture predicts degeneracy. However, no method is known to prove degeneracy in general; even in the  special case of equation \eqref{E.singola} we do not know how to prove that the $S$-integral solutions are not Zariski-dense.

One might try to prove the weaker statement that among the rational points of $\calS_3$ constructed via elliptic fibrations in Section \ref{sec:3pts} only finitely many of them (or maybe only those in a degenerate set) can have integral distances from the three given points.


\end{document}